\input ./style/arxiv-general.cfg
\documentclass[MSNbibl,number,citesort,seceqn,dvips]{arxbj}
\makeatletter
   \@ifpackageloaded{graphicx}{}{\usepackage{graphicx}}
\makeatother
\usepackage{upgreek}

\aid{0}
\volume{21}
\issue{4}
\pubyear{2015}
\firstpage{2190}
\lastpage{2216}
\doi{10.3150/14-BEJ640} 
\docsubty{FLA}

\makeatletter

\newcommand{\rrvert}{\vert}
\newcommand{\rrVert}{\Vert}
\newcommand{\llvert}{\vert}
\newcommand{\llVert}{\Vert}

\newcommand{\bbn}{\mathbb{N}}
\newcommand{\bbr}{\mathbb{R}}
\newcommand{\bbe}{\mathbb{E}}

\newcommand{\cals}{\mathcal{S}}
\newcommand{\calh}{\mathcal{H}}
\newcommand{\calf}{\mathcal{F}}
\newcommand{\cale}{\mathcal{E}}
\newcommand{\calp}{\mathcal{P}}
\newcommand{\calb}{\mathcal{B}}
\newcommand{\calo}{\mathcal{O}}
\newcommand{\calr}{\mathcal{R}}
\newcommand{\calm}{\mathcal{M}}
\newcommand{\cald}{\mathcal{D}}
\newcommand{\Leb}{\operatorname{Leb}}

\newcommand{\La}{{\Lambda}}
\newcommand{\eps}{{\epsilon}}
\newcommand{\ga}{{\gamma}}
\newcommand{\vp}{{\varphi}}
\newcommand{\si}{{\sigma}}
\newcommand{\Om}{{\Omega}}

\newcommand{\dd}{\mathrm{d}}
\newcommand{\ee}{\mathrm{e}}
\newcommand{\ii}{\mathrm{i}}
\newcommand{\bb}{\mathrm{b}}
\newcommand{\cc}{\mathrm{c}}
\newcommand{\cov}{\operatorname{Cov}}

\newtheorem{Theorem}{Theorem}[section]
\newtheorem{Corollary}[Theorem]{Corollary}
\newtheorem{Lemma}[Theorem]{Lemma}
\newproclaim{Definition}[Theorem]{Definition}
\newremark{Remark}[Theorem]{Remark}
\newremark{Example}[Theorem]{Example}
\newcommand{\eqref}[1]{(\ref{#1})}

\renewcommand{\epsilon}{\varepsilon}
\renewcommand{\emptyset}{\varnothing}

\makeatother

\begin{document}
\begin{frontmatter}

\title{Integrability conditions for space--time stochastic integrals:
Theory and applications}
\runtitle{Integrability conditions for space--time integrals}

\begin{aug}
\author[A]{\inits{C.}\fnms{Carsten}~\snm{Chong}\corref{}\thanksref{e1}\ead[label=e1,mark]{carsten.chong@tum.de}\ead[label=u1,url]{www.statistics.ma.tum.de}} \and
\author[A]{\inits{C.}\fnms{Claudia}~\snm{Kl\"uppelberg}\thanksref{e2}\ead[label=e2,mark]{cklu@ma.tum.de}}
\address[A]{Center for Mathematical Sciences, Technische Universit\"at
M\"unchen, Boltzmannstra\ss e 3, 85748 Garching, Germany. \printead{e1,e2};\\ url: \printead*{u1}}
\end{aug}

\received{\smonth{3} \syear{2013}}
\revised{\smonth{3} \syear{2014}}

%
\begin{abstract}
We derive explicit integrability conditions for stochastic integrals
taken over time and space driven by a random measure. Our main tool is
a canonical decomposition of a random measure which extends the results
from the purely temporal case. We show that the characteristics of this
decomposition can be chosen as predictable strict random measures, and
we compute the characteristics of the stochastic integral process. We
apply our conditions to a variety of examples, in particular to ambit
processes, which represent a rich model class.
\end{abstract}

%
\end{frontmatter}

\section{Introduction}\label{s1}

Following It\^o's seminal paper \citep{Ito44}, stochastic integration
theory w.r.t. semimartingales was brought to maturity during the 1970s
and 1980s. One of the fundamental results in this area is the
Bichteler--Dellacherie theorem, which shows the equivalence between the
class of semimartingales and the class of finite $L^0$-random measures.
As a consequence, semimartingales constitute the largest class of
integrators that allow for stochastic integrals of predictable
integrands satisfying the dominated convergence theorem. The natural
analogue to semimartingale integrals in a space--time setting are
integrals of the form
%
\begin{equation}
\label{integral} \int_{\bbr\times E} H(t, x) M(\dd t, \dd x),
\end{equation}
where $E$ is some space and $M$ is an $L^0$-random measure on $\bbr
\times E$. The construction of such integrals is discussed in \citep
{Bichteler83} in its full generality, so the theory is complete from
this point of view.

However, whether $H$ is integrable w.r.t. $M$ or not, depends on whether
%
\begin{equation}
\label{L0crit} \lim_{r\to0} \sup \biggl\{\bbe \biggl[\biggl\llvert
\int S \,\dd M\biggr\rrvert \wedge1 \biggr]\colon|S|\leq|rH|, S\mbox{ is a simple
integrand} \biggr\}=0
\end{equation}
or not, a property which is hard to check. Thus, the aim of this paper
is to characterize \eqref{L0crit} in terms of equivalent conditions,
which can be verified in concrete situations. In the purely temporal
case, this subject is addressed in \citep{Basse13}. The result there
is obtained by using the local semimartingale characteristics
corresponding to a random measure. Our approach parallels this method,
but it turns out that the notion of characteristics in the space--time
setting is much more complex. We will show that, if $M$ has different
times of discontinuity (cf. Definition~\ref{rmclasses} below), we can
associate a characteristic triplet to it consisting of strict random
measures (cf. Definition~\ref{ranmeas}(3)) that are jointly $\si
$-additive in space and time. Moreover, we will determine the
characteristics of stochastic integral processes, which is more
involved than in the temporal case, since a concept is needed to merge
space and time appropriately. Having achieved this step, integrability
conditions in the same fashion as in \citep{Basse13,Rajput89} can be
given for space--time integrals. We will also compare our results to
those of \citep{Rajput89}, \citep{Walsh86} and \citep{Jacod03}.

Applications of our theoretical results will be chosen from the class
of ambit processes
%
\begin{equation}
\label{ambit0} Y(t,x):=\int_{\bbr\times\bbr^d} h(t,s;x,y)\si(s,y) M(\dd s,
\dd y),\qquad  t\in\bbr,x\in\bbr^d,
\end{equation}
which have been suggested for modelling physical space--time phenomena
like turbulence, see, for example, \citep{BN11-2}. In the case, where
$\si=1$ and $M$ is a L\'evy basis (see Remark~\ref{Levybasis}), such
multiparameter integrals have already been investigated by many
authors: for instance, \citep{Cambanis90,Marcus05,Rosinski89} discuss
path properties of the resulting process $Y$, while \citep
{Fasen09,Moser13} address the extremal behaviour of $Y$; mixing
conditions are examined in \citep{Fuchs13}.

As a broad model class, the applications of ambit processes go far
beyond turbulence modelling. For example, \citep{Nagasawa99} describes
the movement of relativistic quantum particles by equations of the form
\eqref{ambit0}. Moreover, solutions to stochastic partial differential
equations driven by random noise are often of the form \eqref{ambit0},
cf. \citep{BN11-2,Walsh86} and Section~\ref{s52}. Furthermore,
stochastic processes like forward contracts in bond and electricity
markets based on a Heath--Jarrow--Morten approach also rely on a spatial
structure, cf. \citep{BN13,Andresen14}. Other applications include
brain imaging \citep{Jensen13} and tumor growth \citep{BN07,Jensen07}.

The concept of an ambit process has also been successfully invoked to
define superpositions of stochastic processes like Ornstein--Uhlenbeck
processes or, more generally, continuous-time ARMA (CARMA) processes.
In these models, only integrals of deterministic integrands w.r.t. L\'
evy bases are involved, so the integration theory of \citep{Rajput89}
is sufficient. Our integrability conditions, however, allow for a
volatility modulation of the noise, which generates a greater model
flexibility. Moreover, in \citep{Behme13} ambit processes have been
used to define superpositions of continuous-time GARCH (COGARCH)
processes. In its simplest case superposition leads to multi-factor
models, economically and statistically necessary extensions of the
one-factor models; cf. \citep{Jacod12}. As we shall see, the
supCOGARCH model again needs the integrability criteria we have
developed since for this model the volatility $\sigma$ and the random
measure $M$ are not independent.

Our paper is organized as follows. Section~\ref{s2} introduces the
notation and gives a summary on the concept of a random measure and its
stochastic integration theory. Section~\ref{s3} derives a canonical
decomposition for random measures as known for semimartingales and
calculates the characteristic triplet of stochastic integral processes.
Section~\ref{s4} presents integrability conditions in terms of the
characteristics from Section~\ref{s3}. Section~\ref{s5} is dedicated
to examples to highlight our results.

\section{Preliminaries}\label{s2}
Let $(\Om,\calf,(\calf_t)_{t\in\bbr},P)$ be a stochastic basis
satisfying the usual assumptions of completeness and right-continuity.
Denote the base space by $\bar\Om:=\Om\times\bbr$ and the optional
(resp. predictable) $\si$-field on $\bar\Om$ by $\calo$ (resp.
$\calp$). Furthermore, fix some Lusin space $E$, equipped with its
Borel $\si$-field $\cale$. Using the abbreviations $\tilde\Om:= \Om
\times\bbr\times E$ and $\tilde\calo:=\calo\otimes\cale$ (resp.
$\tilde\calp:=\calp\otimes\cale$), we call a function $H\dvtx
\tilde\Om\to\bbr$ \textit{optional} (resp. \textit{predictable})
if it
is $\tilde\calo$-measurable (resp. $\tilde\calp$-measurable). We
will often use the symbols $\calo$ and $\calp$ (resp. $\tilde\calo$
and $\tilde\calp$) also for the collection of optional and
predictable functions from $\bar\Om$ (resp. $\tilde\Om$) to $\bbr
$. We refer to \citep{Jacod03}, Chapter I and II,  for all notions not
explicitly explained.

Some further notational conventions: we write $A_t:= A\cap(\Om\times
(-\infty,t])$ for $A\in\calp$, and $\tilde A_t:=\tilde{A}\cap(\Om\times
(-\infty,t]\times E)$ for $\tilde A\in\tilde\calp$. $\calb_\bb
(\bbr^d)$ denotes the collection of bounded Borel sets in $\bbr^d$.
Next, if $\mu$ is a signed measure and $X$ a finite variation process,
we write $|\mu|$ and $|X|$ for the variation of $\mu$ and the
variation process of $X$, respectively. Finally, we equip $L^p=L^p(\Om
,\calf,P)$, $p\in[0,\infty)$, with the topology induced by
\begin{eqnarray*}
\|X\|_p &:=& \bbe\bigl[|X|^p\bigr]^{1/p},\qquad  p\geq1,\\
\|X\|_p& :=& \bbe\bigl[|X|^p\bigr],\qquad  0<p<1, \qquad \|X
\|_0 := \bbe[|X|\wedge1]
\end{eqnarray*}
for $X\in L^p$. Among several definitions of a random measure in the
literature, the following two are the most frequent ones: in essence, a
random measure is either a random variable whose realizations are
measures on some measurable space (e.g., \citep{Jacod03,Kallenberg86})
or it is a $\si$-additive set function with values in the space $L^p$
(e.g., \citep{Bichteler83,Kurtz96,Metivier77,Rajput89,Walsh86}). Our
terminology is as follows:

\begin{Definition}\label{ranmeas}
Let $(\tilde O_k)_{k\in\bbn}$ be a sequence of sets in $\tilde\calp
$ with $\tilde O_k\uparrow\tilde\Om$. Set $\tilde\calp_M:=\bigcup_{k=1}^\infty\tilde\calp|_{\tilde O_k}$, which is the collection of
all sets $A\in\tilde\calp$ such that $A\subseteq\tilde O_k$ for
some $k\in\bbn$.
\begin{enumerate}[(3)]
\item[(1)] An \emph{$L^p$-random measure} on $\bbr\times E$ is a mapping
$M\dvtx \tilde\calp_M \to L^p$ satisfying:
\begin{enumerate}[(a)]
\item[(a)]$M(\emptyset)=0$ a.s.,
\item[(b)] For every sequence $(A_i)_{i\in\bbn}$ of pairwise disjoint sets
in $\tilde\calp_M$ with $\bigcup_{i=1}^\infty A_i\in\tilde\calp
_M$ we have
\[
M \Biggl(\bigcup_{i=1}^\infty
A_i \Biggr)=\sum_{i=1}^\infty
M(A_i) \qquad \mbox{in } L^p.
\]
\item[(c)] For all $A\in\tilde\calp_M$ with $A\subseteq\tilde\Om_t$
for some $t\in\bbr$, the random variable $M(A)$ is $\calf_t$-measurable.
\item[(d)] For all $A\in\tilde\calp_M$, $t\in\bbr$ and $F\in\calf_t$,
we have
\[
M \bigl(A\cap\bigl(F\times(t,\infty)\times E\bigr) \bigr)=1_F M
\bigl(A\cap\bigl(\Om \times(t,\infty)\times E\bigr) \bigr)\qquad  \mbox{a.s.}
\]
\end{enumerate}
\item[(2)] If $p=0$, we only say \emph{random measure}; if $\tilde O_k$ can
be chosen as $\tilde\Om$ for all $k\in\bbn$, $M$ is called a
\emph
{finite} random measure; and finally, if $E$ consists of only one
point, $M$ is called a \emph{null-spatial} random measure.
\item[(3)] A \emph{strict random measure} is a signed transition kernel
$\mu
({\omega},\dd t,\dd x)$ from $(\Om,\calf)$ to\linebreak[4]  $(\bbr\times E, \calb
(\bbr)\otimes\cale)$ with the following properties:
\begin{enumerate}[(a)]
\item[(a)] There is a strictly positive function $V\in\tilde\calp$ such
that $\int_{\bbr\times E} V(t,x) |\mu|(\dd t,\dd x)\in L^1$.
\item[(b)] For $\tilde\calo$-measurable functions $W$ such that $W/V$ is
bounded, the process
\[
W\ast\mu_t:=\int_{(-\infty,t]\times E} W(s,x) \mu(\dd s,\dd x),\qquad
t\in\bbr,
\]
is optional.
\end{enumerate}
\end{enumerate}
\end{Definition}

%
\begin{Remark}
\begin{enumerate}[(2)]
\item[(1)] If we can choose $O_k=\Om\times O^\prime_k$ with $O^\prime
_k\uparrow\bbr\times E$, one popular choice for $(\calf_t)_{t\in
\bbr}$ is the \textit{natural filtration} $(\calf_t^M)_{t\in\bbr}$ of
$M$ which is the smallest filtration satisfying the usual assumptions
such that for all $t\in\bbr$ we have $M(\Om\times B)\in\calf^M_t$
if $B\subseteq((-\infty,t]\times E)\cap O_k^\prime$ with some $k\in
\bbn$.
\item[(2)] If $\mu$ is a positive transition kernel in Definition~\ref
{ranmeas}(3), $\mu$ is an optional $\tilde\calp$-$\si$-finite
random measure in the sense of \citep{Jacod03}, Chapter II, where also
the predictable compensator of a strict random measure is defined.
Obviously, a strict random measure \textit{is} a random measure. For more
details on that, see also \citep{Bichteler83}, Examples 5 and 6.
\end{enumerate}
\end{Remark}

Stochastic integration theory in space--time w.r.t. $L^p$-random
measures is discussed in \citep{Bichteler83}, see also \citep
{Bichteler02}. The special case of $L^2$-integration theory is also
discussed in \citep{Doob53,Walsh86}. Let us recall the details
involved: a \textit{simple integrand} is a function $\tilde\Om\to
\bbr
$ of the form
%
\begin{equation}
\label{simple} S:=\sum_{i=1}^r
a_i1_{A_i},\qquad  r\in\bbn , a_i\in \bbr,
A_i\in\tilde\calp_M,
\end{equation}
for which the stochastic integral w.r.t. $M$ is canonically defined as
%
\begin{equation}
\label{simpleint} \int S \,\dd M := \sum_{i=1}^r
a_i M(A_i).
\end{equation}
Now consider the collection $\cals_M^\uparrow$ of positive functions
$\tilde\Om\to\bbr$ which are the pointwise supremum of simple
integrands and define the \textit{Daniell mean} $\|\cdot\|^\mathrm
{D}_{M,p}\dvtx \bbr^{\tilde\Om} \to[0,\infty]$ by
\begin{itemize}
\item$\|K\|^\mathrm{D}_{M,p} := \sup_{S\in\cals_M,|S|\leq K} \llVert \int S \,\dd M\rrVert _p$, if $K\in\cals_M^\uparrow$, and
\item$\|H\|^\mathrm{D}_{M,p} := \inf_{K\in\cals_M^\uparrow
,|H|\leq K} \|K\|^\mathrm{D}_{M,p}$ for arbitrary functions $H\dvtx
\tilde\Om\to\bbr$.
\end{itemize}
An arbitrary function $H\dvtx \tilde\Om\to\bbr$ is called \textit
{integrable} w.r.t. $M$ if there is a sequence of simple integrands
$(S_n)_{n\in\bbn}$ such that $\|H-S_n\|^\mathrm{D}_{M,p} \to0$ as
$n\to\infty$. Then the \textit{stochastic integral} of $H$ w.r.t. $M$
defined by
%
\begin{equation}
\label{defint} \int H \,\dd M := \int_{\bbr\times E} H(t,x) M(\dd t,\dd x)
:= \lim_{n\to\infty} \int S_n \,\dd M
\end{equation}
exists in $L^p$ and does not depend on the choice of $(S_n)_{n\in\bbn
}$. The collection of integrable functions is denoted by $L^{1,p}(M)$
and can be characterized as follows (\citep{Bichteler02}, Theorems 3.4.10 and
3.2.24):

%
\begin{Theorem}\label{L1p-char} Let $F^{1,p}(M)$ be the collection of functions
$H$ with $\|rH\|^{\mathrm{D}}_{M,p} \to0$ as $r\to0$.\vspace*{2pt} If we identify
two functions coinciding up to a set whose indicator function has
Daniell mean $0$, then
%
\begin{equation}
L^{1,p}(M)=\tilde\calp\cap F^{1,p}(M).
\end{equation}
Moreover, the following dominated convergence theorem holds:
Let $(H_n)_{n\in\bbn}$ be a sequence in $L^{1,p}(M)$ converging
pointwise to some limit $H$. If there exists some function $F\in
F^{1,p}(M)$ with $|H_n|\leq F$ for each $n\in\bbn$, both $H$ and
$H_n$ are integrable with $\|H-H_n\|^{\mathrm{D}}_{M,p}\to0$ as $n\to
\infty$ and
{\renewcommand{\theequation}{DCT}
%
\begin{equation}
\label{DCT} \int H \,\dd M = \lim_{n\to
\infty} \int H_n
\,\dd M \qquad \mbox{in } L^p.
\end{equation}}
\end{Theorem}

Given a predictable function $H\in\tilde\calp$, we can obviously
define a new random measure $H.M$ in the following way:
\setcounter{equation}{4}
\begin{equation}
\label{HM} K\in L^{1,0}(H.M) :\Leftrightarrow KH\in
L^{1,0}(M), \qquad \int K \,\dd(H.M) := \int KH \,\dd M.
\end{equation}
This indeed defines a random measure provided there exists a sequence
$(\tilde O_k)_{k\in\bbn}\subseteq\tilde\calp$ with\linebreak[4]  $\tilde O_k
\uparrow\tilde\Om$ and $1_{\tilde O_k}\in L^{1,0}(H.M)$ for all
$k\in\bbn$. But this construction does not extend the class
$L^{1,0}(M)$ of integrable functions w.r.t. $M$. However, as shown in
\citep{Bichteler83}, Section 3,  $L^{1,p}(M)$ can indeed be extended further
in the following way. Given an $L^p$-random measure $M$, fix some
$\tilde\calp$-measurable function $H$ such that:
%
\begin{equation}
\label{regularity} \mbox{There exists a predictable process } K\dvtx \bar\Om\to\bbr,
K>0, \mbox{ such that } KH\in L^{1,p}(M).
\end{equation}
Now set $\bar O_k:=\{K\geq k^{-1}\}$ for $k\in\bbn$, which obviously
defines predictable sets increasing  to $\bar\Om$, and then $\calp
_{H\cdot M}:=\{A\in\calp\colon A\subseteq\bar O_k \mbox{ for some }
k\in\bbn\}$. Then we define a new null-spatial $L^p$-random measure by
\[
H\cdot M\dvtx \calp_{H\cdot M} \to L^p,\qquad (H\cdot M) (A) :=
\int1_A H \,\dd M.
\]
The following is known from \citep{Bichteler83}, see also \citep{Basse13}, Theorem
A.4:
\begin{enumerate}[(2)]
\item[(1)] If $H\in L^{1,p}(M)$, $H\cdot M$ is a finite $L^p$-random measure
and $\int1 \,\dd(H\cdot M)=\int H \,\dd M$.
\item[(2)] If $G\dvtx \bar\Om\to\bbr$ is a predictable process, we have
$G\in L^{1,p}(H\cdot M)$ if and only if $\|rGH\|_{M,p}\to0$ as $r\to
0$, where for every $\tilde\calp$-measurable function $H$ we set
%
\begin{equation}
\label{eq1} \|H\|_{M,p} := \mathop{\sup_{
F\colon\bar\Om\to\bbr\mathrm{\ predictable},}}_{
|F|\leq1, FH\in L^{1,p}(M)}
\biggl\llVert \int FH \,\dd M\biggr\rrVert _p.
\end{equation}
In this case, we have $\int G\, \dd(H\cdot M) = \int GH\, \dd M$.
\end{enumerate}
Therefore, it is reasonable to extend the set of \textit{integrable}
functions w.r.t. $M$ from $L^{1,p}(M)$ to
%
\begin{equation}
\label{Lpvar} L^p(M) = \bigl\{H\in\tilde\calp\colon H \mbox{ satisfies
\eqref{regularity}  and } \|rH\|_{M,p}\stackrel {r\to0} {
\longrightarrow}0\bigr\}
\end{equation}
by setting
\[
\int H \,\dd M := (H\cdot M) (\bar\Om),\qquad  H\in L^p(M).
\]
We remark that in the null-spatial case $L^{1,0}(M)= L^0(M)$. But in
general, the inclusion $L^{1,p}(M)\subseteq L^p(M)$ is strict, see
\citep{Bichteler83}, Section 3b,  and Example~\ref{LambdaH} below.

Let us also remark that \citep{Dalang99} introduces a stochastic
integral for a Gaussian random measure where the integrands are allowed
to be distribution-valued. It is still an open question whether it is
possible to extend this to the general setting of $L^p$-random
measures, in particular if $p<2$; we do not pursue this direction in
the present paper.

In the sequel we will frequently use the following fact from \citep{Basse13}, Example
3.1: If $M$ is a finite random measure, the process
$(M(\tilde\Om_t))_{t\in\bbr}$ has a c\`adl\`ag modification, which
is then a semimartingale up to infinity w.r.t. to the underlying
filtration (see \citep{Basse13}, Section~2,  for a definition). This
semimartingale will be also be denoted by $M=(M_t)_{t\in\bbr}$.

\section{Predictable characteristics of random measures}\label{s3}
Let us introduce three important subclasses of random measures:

\begin{Definition}\label{rmclasses}
Let $M$ be a random measure where $\tilde O_k=O_k\times E_k$ with $O_k
\uparrow\bar\Om$ and $E_k \uparrow E$. Set $\cale_M:=\bigcup_{k=1}^\infty\cale|_{E_k}$.
\begin{enumerate}[(3)]
\item[(1)]$M$ has \emph{different times of discontinuity} if for all
$k\in
\bbn$ and disjoint sets $U_1, U_2 \in\cale_M$ the semimartingales
$1_{O_k\times U_i}\cdot M$, $i=1,2$, a.s. never jump at the same time.
\item[(2)]$M$ is called \emph{orthogonal} if for all pairs of disjoint sets
$U_1,U_2\in\cale_M$ and $k\in\bbn$ we have $[(1_{O_k\times
U_1}\cdot M)^\cc,(1_{O_k\times U_2}\cdot M)^\cc]=0$.
\item[(3)]$M$ has \emph{no fixed time of discontinuity} if for all $U\in
\cale_M$, $k\in\bbn$ and $t\in\bbr$ we have\linebreak[4]  $\Delta(1_{O_k\times
U}\cdot M)_t=0$ a.s.
\end{enumerate}
\end{Definition}

In the next theorem, we prove a canonical decomposition for random
measures with different times of discontinuity generalizing the results
of \citep{Jacod03} and \citep{Basse13}. Without this extra assumption
on the random measure, only non-explicit results such as \citep{Bichteler83}, Theorem
4.21,  or results for $p\geq2$ as in \citep{Lebedev95b}, Theorem
1,  are known. We also remark that the integrability
conditions in Theorem~\ref{main} will be stated in terms of this decomposition.
Some notation beforehand: we write $\calb_0(\bbr)$ for the collection
of Borel sets on $\bbr$ which are bounded away from $0$. Furthermore,
if $X$ is a semimartingale up to infinity, we write $\mathfrak{B}(X)$
for its first characteristic, $[X]$ for its quadratic variation, $X^\cc
$ for its continuous part (all of them starting at $-\infty$ with
$0$), $\mu^X$ for its jump measure and $\nu^X$ for its predictable
compensator. Finally, if $U\in\cale$, $M|_U$ denotes the random
measure given by $M|_U(A)=M(A\cap(\bar\Om\times U))$ for $A\in
\tilde\calp_M$.

%
\begin{Theorem}\label{candec}
Let $M$ have different times of discontinuity.
\begin{enumerate}[(3)]
\item[(1)] The mappings
\[
B(A):=\mathfrak{B}(1_A\cdot M)_\infty,\qquad  M^\cc(A):=(1_A
\cdot M)^\cc _\infty,\qquad  A\in\tilde\calp_M,
\]
are random measures, the mapping
\[
C(A;B):=\bigl[(1_A\cdot M)^\cc,(1_B\cdot
M)^{\cc}\bigr]_\infty,\qquad  A\in\tilde \calp_M,
\]
is a random bimeasure (i.e., a random measure in both arguments when the
other one is fixed) and
%
\begin{equation}
\label{munu} \hspace*{-15pt}\mu(A,V):=\mu^{1_A\cdot M}(\bbr\times V), \nu (A,V):=
\nu^{1_A\cdot M}(\bbr\times V), \qquad A\in\tilde\calp_M, V\in
\calb_0(\bbr),
\end{equation}
can be extended to random measures on $\tilde\calp_M\otimes\calb
_0(\bbr)$.
Moreover, $(B,C,\nu)$ can be chosen as predictable strict random
(bi-)measures and form the \emph{characteristic triplet} of $M$.
\item[(2)] Let $A\in\tilde\calp_M$ and $\tau$ be a truncation function
(i.e., a bounded function with $\tau(y)=y$ in a neighbourhood of $0$).
Then $1_A(t,x)(y-\tau(y))$ (resp. $1_A(t,x)\tau(y)$) is integrable
w.r.t. $\mu$ (resp. $\mu-\nu$), and we have
%
\begin{eqnarray}
\label{dec} M(A)&=&B(A)+M^{\cc}(A) + \int_{\bbr\times
E \times\bbr}
1_A(t,x) \bigl(y-\tau(y)\bigr) \mu(\dd t,\dd x,\dd y)
\nonumber
\\[-8pt]\\[-8pt]
&&{}+ \int_{\bbr\times E \times\bbr} 1_A(t,x)\tau(y) (\mu-\nu) (\dd t,
\dd x,\dd y),\nonumber
\end{eqnarray}
\item[(3)] There are a positive predictable strict random measure $A({\omega
},\dd t,\dd x)$, a $\tilde\calp$-measurable function $b({\omega
},t,x)$ and a transition kernel $K({\omega},t,x,\dd y)$ from $(\tilde
\Om,\tilde\calp)$ to $(\bbr,\calb(\bbr))$ such that for a.e.
${\omega}\in\Om$
\begin{eqnarray*}
B({\omega},\dd t,\dd x)&=&b({\omega},t,x) A({\omega},\dd t,\dd x), \\
 \nu ({\omega},
\dd t,\dd x, \dd y)&=&K({\omega},t,x,\dd y) A({\omega}, \dd t,\dd x).
\end{eqnarray*}
\end{enumerate}
\end{Theorem}

For the proof of Theorem~\ref{candec}, let us recall the
semimartingale topology of \citep{Emery79} on the space $\cals\calm$
of semimartingales up to infinity, which is induced by
\[
\|X\|_{\mathcal{S}\mathcal{M}} := \sup_{|H|\leq1, H\in\calp} \biggl\llVert \int
_{-\infty}^\infty H_t \,\dd X_t
\biggr\rrVert _0,\qquad  X\in\cals \calm.
\]
The following results are known.

\begin{Lemma}\label{contdec}
\begin{enumerate}[(3)]
\item[(1)] Let $(X^n)_{n\in\bbn}\subseteq\cals\calm$ and $(B^n,C^n,\nu
^n)$ denote the semimartingale characteristics of $X^n$. If $X^n\to0$
in $\mathcal{S}\mathcal{M}$, then each of the following
semimartingale sequences converges to $0$ in $\cals\calm$ as well:
$B^n$, $X^{\mathrm{c},n}$, $C^n$, $[X^n]$, $(y-\tau(y))\ast\mu^n$
and $\tau(y)\ast(\mu^n-\nu^n)$.
\item[(2)] If $W({\omega},t,y)$ is a positive bounded predictable function,
then $W\ast\mu^n\to0$ in probability if and only if $W\ast\nu^n\to
0$ in probability. Morever, $W\ast\nu^n<\infty$ a.s. implies $W\ast\mu^n<\infty$ a.s.
\item[(3)] The collection of predictable finite variation processes is
closed under the semimartingale topology.
\end{enumerate}
\end{Lemma}

For the first part of this lemma, see \citep{Basse13}, Theorem 3.5,  and
\citep{Emery79}, page 276. The second part is taken from \citep{Basse13}, Lemmas
3.1 and 3.3. The third assertion is proved in \citep{Memin80}, Theorem
IV.7.

\begin{pf*}{Proof of Theorem~\ref{candec}}
Let $k\in\bbn$ and consider the set function $(S,U) \mapsto B(S\times
U)$ from the semiring $\calh:=\calp|_{O_k}\times\cale|_{E_k}$ to
$L^0$. Obviously, it is finitely additive in each component: for fixed
$U$, additivity in time holds by the definition of $B$, while for fixed
$S$, additivity in space is due to the assumption of different times of
discontinuity. By a straightforward induction argument this implies
that $B$ is also finitely additive jointly in space and time. Next, let
\[
\calr(\calh)= \Biggl\{ \bigcup_{n=1}^N
C_n\colon N\in\bbn, C_n\in \calh\mbox{ pairwise
disjoint} \Biggr\}
\]
denote the ring generated by $\calh$. Setting $B(\bigcup_{n=1}^N
C_n):=\sum_{n=1}^N B(C_n)$ one obtains a well-defined extension of $B$
to $\calr(\calh)$, which is consistent with the original definition
of $B$ and still finitely additive. Furthermore, since $\calr(\calh)$
contains $O_k\times E_k$, we can further extend $B$ to a measure on
$\si(\calh)=\tilde\calp|_{\tilde O_k}$ using \citep{Kwapien92}, Theorem
B.1.1. We only have to show the implication
%
\begin{equation}
\label{cex} (A_n)_{n\in\bbn}\subseteq\calr(\calh)\qquad  \mbox{with}\quad  \limsup_{n\to\infty} A_n = \emptyset\quad \Longrightarrow\quad \lim
_{n\to
\infty} B(A_n) = 0 \qquad \mbox{in } L^0.
\end{equation}
In fact, under the assumption on the left-hand side of \eqref{cex},
$1_{A_n}\cdot M\to0$ in $\cals\calm$:
\begin{eqnarray*}
\|1_{A_n}\cdot M\|_{\cals\calm} &=& \sup_{|H|\leq1, H\in\calp}
\biggl\llVert \int H\, \dd(1_{A_n}\cdot M)\biggr\rrVert _0 =
\sup_{|H|\leq1, H\in
\calp} \biggl\llVert \int H1_{A_n}\, \dd M\biggr
\rrVert _0
\\
&\leq&\sup_{S\in\cals_M, |S|\leq1_{A_n}} \biggl\llVert \int S \,\dd M\biggr\rrVert
_0 = \|1_{A_n}\|^{\mathrm{D}}_{M,0} \stackrel{n
\to\infty } {\longrightarrow} 0
\end{eqnarray*}
by \eqref{DCT} with $1_{O_k\times E_k}$ as dominating function. Using
Lemma~\ref{contdec}(1), equation \eqref{cex} follows.

This extension still coincides with the definition of $B$ in
Theorem~\ref{candec}: From the construction given in the proof of \citep{Kwapien92}, Theorem
B.1.1, we know that given $A\in\tilde\calp|_{\tilde
O_k}$, there is a sequence of sets $(A_n)_{n\in\bbn}$ in $\calr
(\calh)$ with $\limsup((A\setminus A_n)\cup(A_n\setminus A)) =
\emptyset$ and $B(A_n) \to B(A)$ in $L^0$ as $n\to\infty$. As above
we obtain $1_{A_n}\cdot M \to1_A \cdot M$ in $\cals\calm$, which
implies the assertion. And of course, $B$ is unique and $B(A)$ does not
depend on the choice of $k\in\bbn$ with $A\subseteq O_k$.

Finally, we prove that $B$ corresponds to a predictable strict random
measure. By \citep{Bichteler83}, Theorem 4.10,  it suffices to show that
for $H\in L^{1,0}(B)$ the semimartingale $H\cdot B$ is predictable and
has finite variation on bounded intervals. If $H\in\cals_M$, this
follows from linearity and the fact that the first characteristic of a
semimartingale up to infinity is a predictable finite variation
process. In the general case choose a sequence $(S_n)_{n\in\bbn
}\subseteq\cals_M$ with $S_n \to H$ pointwise and $|S_n|\leq H$ for
all $n\in\bbn$. As $n\to\infty$, we have $S_n \cdot B \to H \cdot
B$ in $\cals\calm$ by \eqref{DCT}. By Lemma~\ref{contdec}(3) we
conclude that also $H\cdot B$ is a predictable finite variation process.

For $C$ we fix one argument and apply the same procedure to the other
argument; for $M^\cc$ we refer to \citep{Bichteler83}, Theorem 4.13. Let
us proceed to $\mu$ and $\nu$, where in both cases we first fix some
$V\in\calb_0(\bbr)$ with $\inf\{|x|\colon x\in V\} \geq\eps>0$
and $\eps<1$. In order to apply the same construction scheme as for
$B$, only the proof of \eqref{cex} is different for $\mu$ and $\nu$.
To this end, let $(A_n)_{n\in\bbn}$ be as on the left-hand side of
\eqref{cex}, that is, $1_{A_n}\cdot M \to0$ in $\cals\calm$. Now
define $\tilde\tau(y) = (y \wedge\eps) \vee(-\eps)$ and choose
$K>1$ such that $|\tilde\tau(y)|\leq K(y^2\wedge1)$ for $|y|\geq
\eps$. Then
\begin{eqnarray*}
\bigl\|\mu(A_n,V)\bigr\|_0&=& \biggl\llVert \frac{1_V(y)}{|\tilde\tau
(y)|}\bigl|
\tilde\tau(y)\bigr|\ast\mu^{1_{A_n}\cdot M}_\infty\biggr\rrVert _0
\leq \eps^{-1}\bigl\llVert 1_V(y)\bigl|\tilde\tau(y)\bigr|\ast
\mu^{1_{A_n}\cdot
M}_\infty\bigr\rrVert _0
\\
&\leq& K\eps^{-1}\bigl\llVert \bigl(y^2\wedge1\bigr)\ast
\mu^{1_{A_n}\cdot M}_\infty \bigr\rrVert _0\leq K
\eps^{-1}\bigl\llVert [1_{A_n}\cdot M]_\infty\bigr
\rrVert _0 \to0,
\end{eqnarray*}
where the last step follows from Lemma~\ref{contdec}(1). Part (2) of
the same lemma yields that also $\nu(A_n,V) \to0$ in $L^0$ as $n\to
\infty$. Consequently, \citep{Bichteler83}, Theorem 4.12,  shows that $\mu
(\cdot,V)$ and $\nu(\cdot, V)$ can be chosen as positive strict
random measures. Observing that $\mu(A,\cdot)$ (resp. $\nu(A,\cdot
)$) is clearly also a positive (and predictable) strict random measure
for given $A\in\tilde\calp_M$, $\mu$ (resp. $\nu$) can be extended
to a positive (and predictable) strict random measure on the product
$\tilde\calp_M\otimes\calb_0(\bbr)$ (see \citep{Rajput89}, Proposition
2.4). Of course, $\nu$ is the predictable compensator of
$\mu$.

The integrability of $1_A(t,x)(y-\tau(y))$ (resp. $1_A(t,x)\tau(y)$)
w.r.t. $\mu$ (resp. $\mu-\nu$) is an obvious consequence of \eqref
{munu} and the corresponding statements in the null-spatial case. The
canonical decomposition of $M$ follows since both sides of \eqref{dec}
are random measures coinciding on $\calh$.

Finally, part (3) of Theorem~\ref{candec} can be proved analogously to
\citep{Jacod03}, Proposition II.2.9.
\end{pf*}

\begin{Remark}
If $M$ is additionally orthogonal, we have $C(A;B)=C(A\cap B;A\cap B)$
for all $A,B\in\tilde\calp_M$. Consequently, we may identify $C$
with $C(A):=[(1_A\cdot M)^\cc]_\infty$ for $A\in\tilde\calp_M$. Of
course, $C$ can then be chosen as a predictable strict random
measure.
\end{Remark}

Next, we calculate the characteristics introduced in Theorem~\ref
{candec} in two concrete situations: first, for the random measure of a
stochastic integral process, and second, for a random measure under an
absolutely continuous change of measure. Although the results in both
cases are comparable with the purely temporal setting, the first task
turns out to be the more difficult one. Moreover, the characteristics
for stochastic integral processes are of particular importance for our
integrability conditions in Section~\ref{s4}.

Beforehand, we need some bimeasure theory: it is well known that
bimeasures cannot be extended to measures on the product $\si$-field
in general and that integration theory w.r.t. bimeasures differs from
integration theory w.r.t. measures. Following \citep{Chang83}, let two
measurable spaces $(\Om_i,\calf_i)$, $i=1, 2$, and a bimeasure $\beta
\colon\calf_1\times\calf_2\to\bbr$ be given. We call a pair
$(f_1,f_2)$ of $\calf_i$-measurable functions $f_i$, $i=1, 2$, \textit
{strictly $\beta$-integrable} if
\begin{enumerate}[(3)]
\item[(1)]$f_1$ (resp. $f_2$) is integrable w.r.t. $\beta(\cdot;B)$ for
all $B\in\calf_2$ (resp. $\beta(A;\cdot)$ for all $A\in\calf_1$),
\item[(2)]$f_2$ is integrable w.r.t. the measure $B\mapsto\int_{\Om_1}
f_1({\omega}_1) \beta(\dd{\omega}_1;B)$ and $f_1$ is integrable
w.r.t. the measure $A \mapsto\int_{\Om_2} f_2({\omega}_2) \beta
(A;\dd{\omega}_2)$,
\item[(3)] for all $A\in\calf_1$ and $B\in\calf_2$, the following
integrals are equal:
%
\begin{equation}
\label{betaint} \int_{A} f_1({
\omega}_1) \biggl(\int_{B} f_2({
\omega}_2) \beta(\dd{\omega}_1;\dd{\omega}_2)
\biggr) = \int_{B} f_2({\omega}_2)
\biggl(\int_{A} f_1({\omega}_1)
\beta(\dd{\omega}_1;\dd{\omega}_2) \biggr).
\end{equation}
\end{enumerate}
The \textit{strict $\beta$-integral} of $(f_1;f_2)$ on $(A;B)$, denoted
by $\int_{(A;B)} (f_1;f_2) \,\dd\beta$, is then defined as the common
value \eqref{betaint}.

The next theorem determines the characteristics of stochastic integral
processes, which is \citep{Jacod03}, Proposition IX.5.3,  in the null-spatial case.

%
\begin{Theorem}\label{charint}
Let $M$ be a random measure with different times of discontinuity and
$H\in\tilde\calp$ satisfy \eqref{regularity} with some $K>0$. Then
the null-spatial random measure $H\cdot M$ has characteristics
$(B^{H\cdot M}, C^{H\cdot M}, \nu^{H\cdot M})$ given by
%
\begin{eqnarray}
B^{H\cdot M}(A)&=&(H\cdot B) (A)\nonumber\\[-8pt]\\[-8pt]
&&{}+\int_{\bbr\times E\times\bbr
}1_A(t)
\bigl[\tau\bigl(H(t,x)y\bigr)-H(t,x)\tau(y)\bigr] \nu(\dd t,\dd x,\dd y),\nonumber
\\
C^{H\cdot M}(A)&=&\int_\bbr K^{-2}_t
\,\dd \biggl(\int_{(A_t\times
E;A_t\times E)} (HK;HK) \,\dd C \biggr),
\\
W(t,y)\ast\nu^{H\cdot M}&=&W\bigl(t,H(t,x)y\bigr)\ast\nu
\end{eqnarray}
for all $A\in\calp_{H\cdot M}$ and $\calp\otimes\calb(\bbr
)$-measurable functions $W$ such that $W(t,y)\ast\nu^{H\cdot M}$ exists.

Moreover, if in addition $M$ is orthogonal, then
%
\begin{equation}
C^{H\cdot M}(\dd t) = \int_E H^2(t,x) C(
\dd t,\dd x).
\end{equation}
\end{Theorem}

\begin{pf}
The second part of this theorem is clear as soon as we
have proved
the first part. Since characteristics are defined locally, we may
assume that $H\in L^{1,0}(M)$.
We first consider the continuous part $C^{H\cdot M}$: to this end, let
$(H_n)_{n\in\bbn}$ be a sequence of simple integrands with $|H_n|\leq
|H|$ for all $n\in\bbn$ and $H_n\to H$ pointwise. Since for simple
integrands the claim follows directly from the definition of $C$ and
the bimeasure integral, we would like to use the \eqref{DCT} and Lemma~\ref{contdec}(1) on the one hand and the dominated convergence theorem
for bimeasure integrals (see \citep{Chang83}, Corollary 2.9) on the other
hand to obtain the result. In order to do so, we only have to show that
$(H;H)$ is strictly $C$-integrable, which means by the symmetry of $C$
the following two points: first, that $H$ is integrable w.r.t. the
measure $A\mapsto C(A;B)=[(1_A\cdot M)^{\cc}, (1_B\cdot M)^{\cc
}]_\infty$ for all $B\in\tilde\calp_{M}$, and second, that $H$ is
integrable w.r.t. the measure $A\mapsto\int H(t,x) \,\dd C(A;\dd t, \dd
x)=[(1_A\cdot M)^{\cc}, (H\cdot M)^{\cc}]_\infty$.

Let $G$ be $1_B$ or $H$. From \citep{Lebedev95b}, Theorem~2 and its
Corollary, we know that there exists a probability measure $Q$
equivalent to $P$ such that $M$ is an $L^2(Q)$-random measure with
$G,H\in L^{1,2}(M;Q)$. Since the bounded sets in $L^0(P)$ are exactly
the bounded sets in $L^0(Q)$, convergence in $\|\cdot\|^\mathrm
{D}_{M,0;P}$ is equivalent to convergence in $\|\cdot\|^\mathrm
{D}_{M,0;Q}$. Similarly, stochastic integrals and predictable quadratic
covariation remain unchanged under $Q$ (cf. \citep{Bichteler02}, Proposition~3.6.20,  and \citep{Jacod03}, Theorem III.3.13).
Consequently, if we write $\ga(A):=[1_A\cdot M^\cc, G\cdot M^\cc
]_\infty$ for $A\in\tilde\calp_{M}$, it suffices to show that
\[
\sup_{S\in\cals_{M}, |S|\leq|rH|} \biggl\llVert \int S \,\dd\ga\biggr\rrVert
_{L^0(Q)} = \sup_{S\in\cals_{M}, |S|\leq|rH|} \bigl\|\bigl[ (S\cdot
M)^{\cc},(G\cdot M)^{\cc}\bigr]_\infty
\bigr\|_{L^0(Q)}\to0 \qquad \mbox{as } r\to0.
\]
Indeed, using Fefferman's inequality (cf. \citep{Bichteler02}, Theorem
4.2.7), we can find a constant $R>0$, which only depends
on $G$, such that
\begin{eqnarray*}
&& \sup_{S\in\cals_{M}, |S|\leq|rH|} \bigl\|\bigl[ (S\cdot M)^{\cc
},(G\cdot
M)^{\cc}\bigr]_\infty \bigr\|_{L^0(Q)}\\
&&\quad \leq R\sup
_{S\in\cals_M,
|S|\leq|rH|} \bbe_Q \bigl[\bigl[ (S\cdot
M)^{\cc}\bigr]_\infty \bigr]^{1/2}
\\
&&\quad = R \sup_{S\in\cals_{M}, |S|\leq|rH|} \bigl\|(S\cdot M)^{\cc}_\infty\bigr\|
_{L^2(Q)}=R\|rH\|_{M^\cc,2;Q}^\mathrm{D}\to0
\end{eqnarray*}
as $r\to0$, which finishes the proof for $C^{H\cdot M}$.

For $B^{H\cdot M}$ and $\nu^{H\cdot M}$, we first take some $D\in
\calp\otimes\calb_0(\bbr)$ and claim that
%
\begin{equation}
\label{eq2} 1_D(s,y)\ast\mu^{H\cdot M} = 1_D
\bigl(s,H(s,x)y\bigr)\ast\mu .
\end{equation}
This identity immediately extends to finite linear combinations of such
indicators and thus, by \eqref{DCT}, also to all functions $W({\omega
},t,y)$ for which $W\ast\mu^{H\cdot M}$ exists. By the definition of
the predictable compensator, this statement also passes to the case
where $\mu$ is replaced by $\nu$.

In order to prove \eqref{eq2}, first observe that the jump process of
the semimartingale $H\cdot M$ up to infinity is given by $\Delta
(H\cdot M)_t=(H\cdot M)(\Om\times\{t\}\times E)$. Furthermore, we can
assume that $D$ does not contain any points in $\bar\Om\times\{0\}$.
Hence, in the case where $H=1_A$ with $A\in\tilde\calp_M$, we have
for all $t\in\bbr$
\[
1_D(s,y)\ast\mu^{H\cdot M}_t =
1_D(s,y)\ast\mu^{1_A\cdot
M}_t=1_D(s,y)1_A(s,x)
\ast\mu_t = 1_D\bigl(s,1_A(s,x)y\bigr)\ast
\mu_t.
\]
Now a similar calculation yields that \eqref{eq2} remains true for all
functions $H\in\cals_M$. Finally, let $H\in L^{1,0}(M)$. By
decomposing $H=H^+-H^-$ into its positive and negative part, we may
assume that $H\geq0$ and choose a sequence $(H_n)_{n\in\bbn}$ of
simple functions with $H_n\uparrow H$ as $n\to\infty$. As we have
already seen in the proof of Theorem~\ref{candec}, we have
$1_D(s,y)\ast\mu^{H_n\cdot M}\to1_D(s,y)\ast\mu^{H\cdot M}$ in
$\cals\calm$. On the other hand, if $D$ is of the form $R\times
(a,b]$ with $R\in\calp$ and $(a,b]\subseteq(0,\infty)$ or of the
form $R \times[a,b)$ with $[a,b)\subseteq(-\infty,0)$, then
$1_D({\omega},s,H_n({\omega},s,x)y) \to1_D({\omega},s,H({\omega
},s,x)y)$ as $n\to\infty$ for every $({\omega},s,x,y)\in\tilde\Om
\times\bbr$, which shows that \eqref{eq2} holds up to
indistinguishability. For general $D$, use Dynkin's $\pi$-${\lambda
}$-lemma (\citep{Billingsley95}, Theorem 3.2).

Finally, we compute $B^{H\cdot M}$. The results up to now yield that
for all $t\in\bbr$,
\begin{eqnarray*}
(H\cdot M)_t-\bigl(y-\tau(y)\bigr)\ast\mu^{H\cdot M}_t
&=& (H\cdot B)_t + \bigl(H\cdot M^{\cc}\bigr)_t+H(s,x)
\bigl(y-\tau(y)\bigr)\ast\mu_t
\\
&&{}+H(s,x)\tau(y)\ast(\mu-\nu)_t-\bigl[H(s,x)y-\tau\bigl(H(s,x)y
\bigr)\bigr]\ast\mu_t.
\end{eqnarray*}
By definition, $B^{H\cdot M}$ is the finite variation part in the
canonical decomposition of this special semimartingale, which exactly
equals $H\cdot B+[\tau(H(t,x)y)-H(t,x)\tau(y)]\ast\nu$.
\end{pf}

Finally, we show a Girsanov-type theorem comparable to \citep{Jacod03}, Theorem
III.3.24,  for semimartingales. First, let us introduce some
notation. We consider another probability measure $P^\prime$ on $(\Om
,\calf,(\calf_t)_{t\in\bbr})$ such that $P^\prime_t:=P^\prime
|_{\calf_t}$ is absolutely continuous w.r.t. $P_t:=P|_{\calf_t}$ for
all $t\in\bbr$. Then denote by $Z$ the unique $P$-martingale such
that $Z\geq0$ identically and $Z_t$ is a version of the Radon--Nikodym
derivative $\dd P^\prime_t/\dd P_t$ for all $t\in\bbr$, cf. \citep{Jacod03}, Theorem III.3.4.

Now let $M$ be a random measure with different times of discontinuity
under the probability measure $P$ with characteristics $(B,C,\nu)$
w.r.t. the truncation function $\tau$. We modify the sequence $(\tilde
O_k)_{k\in\bbn}$ of Definition~\ref{ranmeas}(1) by setting
$\tilde O^\prime_k := \tilde O_k \cap(\Om\times(-k,k]\times E)$ for
$k\in\bbn$ and $\tilde\calp_M^\prime:=\bigcup_{k=1}^\infty\tilde
\calp|_{\tilde O_k^\prime}$. Next, we denote the jump measure of $M$
by $\mu$ and set $M^P_\mu(W):=\bbe_P[W\ast\mu_\infty]$ for all
non-negative $\calf\otimes\calb(\bbr)\otimes\cale\otimes\calb
(\bbr)$-measurable functions $W$. Furthermore, for every such $W$,
there exists an $M^P_\mu$-a.e. unique $\tilde\calp\otimes\calb
(\bbr)$-measurable function $M^P_\mu(W|\tilde\calp\otimes\calb
(\bbr))$ such that
\[
M^P_\mu(WU)=M^P_\mu
\bigl(M^P_\mu\bigl(W|\tilde\calp\otimes\calb(\bbr)\bigr)
U\bigr)\qquad  \mbox{for all $\tilde\calp\otimes\calb(\bbr)$-measurable $U\geq 0$}.
\]
Finally, we set
\begin{eqnarray*}
Y(t,x,y)&:=& M^P_\mu\bigl(Z/Z_{-}1_{\{Z_{-}>0\}}|
\tilde\calp \otimes\calb(\bbr)\bigr) (t,x,y),\qquad  t\in\bbr,x\in E, y\in\bbr,
\\
C^Z(A)&:=&\bigl[\bigl(Z_{-}^{-1}\cdot Z
\bigr)^\cc,(1_A\cdot M)^\cc
\bigr]_\infty,\qquad  A\in \tilde\calp^\prime_M.
\end{eqnarray*}
In the last line, the stochastic integral process $Z_{-}^{-1}\cdot Z$
is meant to start at $t_0$, where $t_0\in\bbr$ is chosen such that
$(1_A\cdot M)^\cc=0$ on $(-\infty,t_0]$. Then $C^Z(A)$ is
well defined by \citep{Jacod03}, Proposition~III.3.5a,  and does not depend on
the choice of $t_0$. Moreover, as in Theorem~\ref{candec}, one shows
that $C^Z$ can be chosen as a positive predictable strict random measure.

The following theorem extends \citep{Jacod03}, Theorem III.3.24,  to the
space--time framework.

%
\begin{Theorem}
Under $P^\prime$, $M$ is also a random measure with different times of
discontinuity (w.r.t. $(\tilde O_k^\prime)_{k\in\bbn}$). Its
$P^\prime$-characteristics $(B^\prime,C^\prime,\nu^\prime)$ w.r.t.
$\tau$ are versions of
\begin{eqnarray*}
B^\prime(\dd t,\dd x)&:=&B(\dd t,\dd x) + C^Z(\dd t,\dd x) +
\tau (y) \bigl(Y(t,x,y)-1\bigr) \nu(\dd t,\dd x,\dd y),
\\
C^\prime(\dd t,\dd x)&:=&C(\dd t,\dd x),
\\
\nu^\prime(\dd t,\dd x,\dd y)&:=& Y(t,x,y) \nu(\dd t,\dd x,\dd y).
\end{eqnarray*}
\end{Theorem}

\begin{pf}
Since each set in $\tilde\calp_M^\prime$ is $\calf_t$-measurable
for some $t\in\bbr$, properties (a), (b) and (d) of Definition~\ref
{ranmeas}(1) still hold under $P^\prime$. Since (c) does not depend on
the underlying probability measure, $M$ is also a random measure under
$\tilde P$. To show that $M$ still has different times of discontinuity
under $P^\prime$, it suffices to notice the following: using the
notation of Definition~\ref{rmclasses}, the event that $1_{O_k\times
U_1}\cdot M$ and $1_{O_k\times U_2}\cdot M$ have a common jump in $\bbr
$ is the union over $n\in\bbn$ of the events that they have a common
jump in $(-\infty,n]$. Since these latter events are $\calf
_n$-measurable, their $P^\prime$-probability is $0$, as desired.
Finally, the characteristics under $P^\prime$ can be derived, up to
obvious changes, exactly as in \citep{Jacod03}, Theorem III.3.24.
\end{pf}

\section{An integrability criterion}\label{s4}

The canonical decomposition of $M$ in Theorem~\ref{candec} together
with Theorem~\ref{charint} enables us to reformulate \eqref{Lpvar} in
terms of conditions only depending on the characteristics of $M$. This
result extends the null-spatial case as found in \citep{Jacod03}, Theorem
III.6.30, \citep{Cherny05}, Theorem 4.5, \citep{Basse13}, Theorem
3.2,  or \citep{Kwapien92}, Theorem 9.4.1. It also generalizes
the results of \citep{Rajput89}, Theorem 2.7,  to predictable integrands
and also to random measures which are not necessarily L\'evy bases. Our
proof mimics the approach in \citep{Basse13}, Theorem 3.2,  and takes care
of the additional spatial structure.

\begin{Theorem}\label{main}
Let $M$ be a random measure with different times of discontinuity whose
characteristics w.r.t. some truncation function $\tau$ are given by
Theorem~\ref{candec}. Furthermore, let $H\in\tilde\calp$ satisfy
\eqref{regularity}. Then $H\in L^0(M)$ if and only if each of the
following conditions is satisfied a.s.:
%
\begin{eqnarray}\label{cond1}
\hspace*{-15pt}\int_{\bbr\times E} \biggl\llvert H(t,x)b(t,x) + \int
_{\bbr} \bigl[\tau\bigl(H(t,x)y\bigr)-H(t,x)\tau(y) \bigr]
K(t,x,\dd y)\biggr\rrvert A(\dd t,\dd x) &<& \infty,
\\
\label
{cond2}\int_\bbr K^{-2}_t \,\dd \biggl(\int
_{((-\infty
,t]\times
E;(-\infty,t]\times E)} (HK;HK) \,\dd C \biggr)&<&\infty,
\\
\label{cond3}
\int_{\bbr\times E} \int_{\bbr} \bigl(1\wedge
\bigl(H(t,x)y\bigr)^2 \bigr) K(t,x,\dd y) A(\dd t,\dd x)& < &\infty.
\end{eqnarray}
If $M$ is additionally orthogonal, the spaces $L^0(M)$ and $L^{1,0}(M)$
are equal and condition \eqref{cond2} is equivalent to
%
\begin{equation}
\label{cond2var} \int_{\bbr\times E} H^2(t,x) C(\dd t,\dd x)
<\infty.
\end{equation}
\end{Theorem}

The following lemma is a straightforward extension of \citep{Rajput89}, Lemma~2.8. We omit its proof.

\begin{Lemma}\label{lem4}
For $t\in\bbr$, $x\in E$ and $a\in\bbr$ define
\[
U(t,x,a):= \biggl\llvert ab(t,x) + \int_\bbr \bigl(\tau(ay)
- a\tau(y) \bigr) K(t,x,\dd y)\biggr\rrvert , \qquad \tilde U(t,x,a):=\sup
_{-1\leq c \leq1} U(t,x,ca).
\]
Then there exists a constant $\kappa>0$ such that
\[
\tilde U(t,x,a)\leq U(t,x,a)+\kappa\int_\bbr \bigl(1
\wedge(ay)^2 \bigr) K(t,x,\dd y).
\]
\end{Lemma}

\begin{pf*}{Proof of Theorem~\ref{main}} We first prove that $H\in L^0(M)$
implies \eqref{cond1}--\eqref{cond3}. Since $H\cdot M$ is a
semimartingale up to infinity, $B^{H\cdot M}(\bbr)$ and $C^{H\cdot
M}(\bbr)$ exist. Thus, Theorem~\ref{charint} gives the first two
conditions. For the last condition observe that $(1\wedge y^2)\ast\nu
^{H\cdot M}_\infty<\infty$ a.s. is equivalent to $(1\wedge y^2)\ast
\mu^{H\cdot M}_\infty<\infty$ a.s., which obviously holds since
$H\cdot M$ is a semimartingale up to infinity. This completes the first
direction of the proof.

For the converse statement, we define $\cald:=\{G\in\calp\colon
|G|\leq1, GH\in L^{1,0}(M)\}$.
By \eqref{Lpvar}
we have to show that the set $\{\int GH \,\dd M\colon G\in\cald\}$ is
bounded in $L^0$ (i.e., bounded in probability) whenever $H$ satisfies
\eqref{cond1}--\eqref{cond3}. By Theorem~\ref{charint},
\[
\int GH \,\dd M = \int GH \,\dd M^{\cc} + \tau(GHy)\ast(\mu-\nu
)_\infty+ \bigl(GHy-\tau(GHy)\bigr)\ast\mu_\infty+
B^{GH\cdot M}(\bbr).
\]
We consider each part of this formula separately and show that each of
the sets
%
\begin{eqnarray}
\label{eq3}\bigl\{B^{GH\cdot M}(\bbr)\colon G&\in&\cald\bigr\},
\\
\label{eq4}\biggl\{\int GH \,\dd M^{\cc}\colon G&\in&\cald\biggr
\},
\\
\label{eq5}\bigl\{\tau(GHy)\ast(\mu-\nu)_\infty\colon G&\in&\cald\bigr
\},
\\
\label{eq6}\bigl\{\bigl(GHy -\tau(GHy)\bigr)\ast\mu_\infty\colon G&\in&\cald\bigr
\}
\end{eqnarray}
is bounded in probability.

If $G\in\cald$ and $\kappa>0$ denotes the constant in Lemma~\ref
{lem4}, \eqref{cond1} and \eqref{cond3} imply
\begin{eqnarray*}
&& \int_{\bbr\times E} U\bigl(t,x,G_t H(t,x)\bigr) A(\dd
t,\dd x) \\
&&\quad \leq\int_{\bbr\times E} \tilde U\bigl(t,x,G_t
H(t,x)\bigr) A(\dd t,\dd x)
 \leq \int_{\bbr\times E} \tilde U\bigl(t,x,H(t,x)\bigr) A(\dd t,\dd
x)
\\
&&\quad \leq \int_{\bbr\times E} U\bigl(t,x,H(t,x)\bigr) A(\dd t, \dd x)+
\kappa\int_{\bbr\times E} \int_\bbr \bigl(1\wedge
\bigl(H(t,x)y\bigr) \bigr) K(t,x,\dd y) A(\dd t, \dd x) < \infty
\end{eqnarray*}
a.s., which shows that \eqref{eq3} is bounded in probability.

Next, consider \eqref{eq4} and fix some $G\in\cald$ for a moment.
Using Lenglart's inequality \citep{Jacod03}, Lemma I.3.30a, we have
for all $\eps$, $\eta>0$
\begin{eqnarray*}
&& P \biggl[\biggl\llvert \int GH \,\dd M^{\cc}\biggr\rrvert\geq\eps
\biggr]\\
&&\quad   \leq P \Bigl[\sup_{t\in\bbr}\bigl|\bigl(GH\cdot M^{\cc}
\bigr) (\bar\Om _t)\bigr|\geq\eps \Bigr]= P \Bigl[\sup_{t\in\bbr}\bigl|\bigl(GH\cdot M^{\cc}\bigr) (
\bar\Om_t)\bigr|^2\geq \eps^2 \Bigr]
\\
&&\quad  \leq
\frac{\eta}{\eps^2}+P \bigl[\bigl[GH\cdot M^{\cc
}\bigr]_\infty\geq
\eta \bigr]= \frac{\eta}{\eps^2}+ P \bigl[G^2K^{-2}\cdot\bigl[KH\cdot
M^\cc \bigr]_\infty\geq\eta \bigr]
\\
&&\quad \leq\frac{\eta}{\eps^2} + P
\bigl[K^{-2}\cdot\bigl[KH\cdot M^\cc\bigr]_\infty\geq
\eta \bigr].
\end{eqnarray*}
Now \eqref{cond2} allows us to make the quantity on the left-hand side
arbitrarily small, independently of $G\in\cald$, by first choosing
$\eta>0$ and then $\eps>0$ large enough.

For \eqref{eq5}, we use the abbreviation $W(t,x,y)=\tau(G_t
H(t,x)y)$. Lenglart's inequality again yields\vspace*{-1pt}
%
\begin{equation}
\label{eq15} P\bigl[\bigl|W\ast(\mu-\nu)_\infty\bigr|\geq\eps\bigr] \leq P \Bigl[
\sup_{t\in\bbr
} \bigl|W\ast(\mu-\nu)_t\bigr|^2\geq
\eps^2 \Bigr]\leq\frac{\eta}{\eps
^2} + P \bigl[\bigl\langle W\ast(\mu-
\nu)\bigr\rangle_\infty\geq\eta \bigr]
\end{equation}
for every $\eps$, $\eta>0$. Furthermore, by Theorem~\ref{charint}
and \citep{Jacod03}, Proposition II.2.17,  we have
\[
\bigl\langle W\ast(\mu-\nu)\bigr\rangle_\infty= \bigl\langle\tau(y)\ast
\bigl(\mu ^{GH\cdot M}-\nu^{GH \cdot M}\bigr)\bigr\rangle_\infty\leq
\tau(y)^2\ast\nu _\infty,
\]
which is finite by \eqref{cond3} yielding the boundedness of \eqref{eq5}.

Next, choose $r$, $\eps>0$ such that $f(y):=r|y|1_{\{|y|>\eps\}}$
satisfies $|y-\tau(y)|\leq f(y)$ for all $y\in\bbr$. Obviously, $f$
is symmetric and increasing on $\bbr_+$ so that
\[
\bigl\llvert \bigl(GHy-\tau(GHy) \bigr)\ast\mu_\infty\bigr\rrvert \leq
f(GHy)\ast\mu_\infty\leq f(Hy)\ast\mu_\infty.
\]
Now the third condition and Lemma~\ref{contdec}(2) imply that
\[
\sum_{t\in\bbr} \bigl(1\wedge\eps^2\bigr)
1_{\{|\Delta(H\cdot M)_t|>\eps\}
}\leq\bigl(1\wedge y^2\bigr)\ast\mu^{H\cdot M}_\infty=
\bigl(1\wedge \bigl(H(t,x)y\bigr)^2 \bigr)\ast\mu_\infty<
\infty
\]
a.s. such that $\{|\Delta(H\cdot M)_t|>\eps\}$ only happens for
finitely many time points. Hence,
\[
f(Hy)\ast\mu_\infty=f(y)\ast\mu^{H\cdot M}_\infty=r\sum
_{t\in
\bbr}\bigl|\Delta(H\cdot M)_t\bigr|
1_{\{|\Delta(H\cdot M)_t|>\eps\}}<\infty
\]
a.s., which implies that the set in \eqref{eq6} is also bounded in probability.

Finally, in the case where $M$ is also orthogonal, we show that \eqref
{cond1}, \eqref{cond2var} and \eqref{cond3} imply $H\in L^{1,0}(M)$.
By Theorem~\ref{L1p-char} and the fact that for predictable functions $H$
\[
\|H\|_{M,0}^{\mathrm{D}} = \sup_{S\in\cals_M, |S|\leq|H|} \biggl
\llVert \int S \,\dd M\biggr\rrVert _0 = \sup_{G\in\tilde\calp, |G|\leq1, GH\in
L^{1,0}(M)}
\biggl\llVert \int GH \,\dd M\biggr\rrVert _0,
\]
we have to show that the set $\{\int GH \,\dd M\colon G\in\cald^\prime
\}$ is bounded in $L^0$, where $\cald^\prime$ consists of all
functions $G\in\tilde\calp$ with $|G|\leq1$ and $GH\in L^{1,0}(M)$.
Obviously, the previously considered set $\cald$ is a subset of $\cald
^\prime$. Intending to verify \eqref{eq3}--\eqref{eq6} with $G$ taken
from $\cald^\prime$, we observe that all calculations remain valid
except those for \eqref{eq4}. For \eqref{eq4} we argument as follows:
for all $\eps$, $\eta>0$, Lenglart's inequality implies
\begin{eqnarray*}
&& P \biggl[\biggl\llvert \int GH \,\dd M^{\cc}\biggr\rrvert \geq\eps
\biggr] \\
&&\quad \leq P \Bigl[\sup_{t\in\bbr}\bigl|(GH\cdot M)^{\cc}(\bar
\Om _t)\bigr|^2\geq\eps^2 \Bigr]\leq \frac{\eta}{\eps^2}+P \bigl[\bigl[ (GH\cdot M)^{\cc}
\bigr]_\infty\geq \eta \bigr]
\\
&&\quad =\frac{\eta}{\eps^2}+ P \biggl[\int
_{\bbr\times E} G^2(t,x) H^2(t,x) C(\dd t, \dd
x)\geq\eta \biggr]
\\
&&\quad \leq \frac{\eta}{\eps^2} + P \biggl[\int_{\bbr\times E}
H^2(t,x) C(\dd t, \dd x)\geq\eta \biggr].
\end{eqnarray*}
This finishes the proof of Theorem~\ref{main}.
\end{pf*}

The remaining part of this section illustrates Theorem~\ref{main} by a
series of remarks, examples and useful extensions.

\begin{Remark}\label{h0}
If $M$ has summable jumps, which means that each of the semimartingales\linebreak[4]
$(M(\tilde\Om_t \cap\tilde O_k))_{t\in\bbr}$, $k\in\bbn$, has
summable jumps over finite intervals, it is often convenient to
construct the characteristics w.r.t. $\tau=0$, which is not a proper
truncation function. Then one would like to use $\tau=0$ in \eqref
{cond1} and replace \eqref{cond3} by
%
\begin{equation}
\label{3cond-var} \int_{\bbr\times E} \int_\bbr
\bigl(1\wedge \bigl|H(t,x)y\bigr| \bigr) K(t,x,\dd y) A(\dd t,\dd x) < \infty.
\end{equation}
We show that \eqref{cond1} with $\tau=0$, \eqref{cond2} and \eqref
{3cond-var} are together sufficient conditions for $H\in L^0(M)$.
First, note that we can choose $\kappa=0$ in Lemma~\ref{lem4}(2)
since $\tau$ is identical $0$ and therefore $\tilde U = U$. So the
calculations done for \eqref{eq3} remain valid. Moreover, \eqref{eq4}
does not depend on $\tau$ and the boundedness of \eqref{eq5} becomes
trivial. For \eqref{eq6} observe that
%
\begin{equation}
\label{eq7} |GHy|\ast\mu_\infty\leq|Hy|\ast\mu _\infty= |y|
\ast\mu^{H\cdot M}_\infty= |y|1_{\{|y|\leq1\}}\ast
\mu^{H\cdot
M}_\infty+ |y|1_{\{|y|>1\}}\ast\mu^{H\cdot M}_\infty.
\end{equation}
Now \eqref{3cond-var} implies by Lemma~\ref{contdec}(2) that a.s.,
\[
|y|1_{\{|y|\leq1\}}\ast\mu^{H\cdot M}_\infty+ 1_{\{|y|>1\}}\ast
\mu ^{H\cdot M}_\infty< \infty.
\]
As a result, on the right-hand side of \eqref{eq7}, the first summand
converges a.s. and the second one is in fact just a finite sum a.s.

The converse statement is not true, already in the null-spatial case:
let $(N_t)_{t\geq0}$ be\vspace*{2pt} a standard Poisson process and $\tilde
N_t=N_t-t$, $t\geq0$, its compensation. Set $H_t:=(1+t)^{-1}$ for
$t\geq0$. Then $H\in L^0(\tilde N)$ as one can see from \eqref
{cond1}--\eqref{cond3} with the proper truncation function $\tau
(y)=y1_{\{|y|<1\}}$; but $\int_0^\infty H_t \,\dd t=\infty$ violating
both \eqref{cond1} with $\tau=0$ and \eqref{3cond-var}.

However, if $M$ is a positive (or negative) random measure, that is,
$M(A)$ is a positive (or negative) random variable for all $A\in\tilde
\calp_M$, then $C=0$ necessarily and \eqref{cond1} with $\tau=0$ and
\eqref{3cond-var} also become necessary conditions for $H\in
L^0(M)=L^{1,0}(M)$; cf. \citep{Bichteler83}, Example 5, page~7, and Theorem
4.12.
\end{Remark}

Next, we compare our results and techniques to the standard literature.

\begin{Remark}[(L\'evy bases \citep{Rajput89})] \label{Levybasis} L\'
evy bases
are originally called infinitely divisible independently scattered
random measures in \citep{Rajput89}. They are the space--time analogues
of processes with independent increments and have attracted interest in
several applications in the last few years, see Section~\ref{s5} for
some examples.
The precise definition is as follows: Assume that we have $\tilde
O_k=\Om\times O^\prime_k$ in the notation of Definition~\ref
{ranmeas}, where $(O^\prime_k)_{k\in\bbn}$ is a sequence increasing
to $\bbr\times E$. Set $\cals:=\bigcup_{k=1}^\infty\calb(\bbr
^{1+d})|_{O^\prime_k}$. Then a \textit{L\'evy basis} $\La$ is a random
measure on $\bbr\times E$ with the following additional properties:
\begin{enumerate}[(2)]
\item[(1)] If $(A_n)_{n\in\bbn}$ is a sequence of pairwise disjoint sets
in $\cals$, then $(\La(\Om\times A_n))_{n\in\bbn}$ are independent
random variables.
\item[(2)] For all $A\in\cals$, $\La(\Om\times A)$ has an infinitely
divisible distribution.
\end{enumerate}

Note that we have altered the original definition of \citep{Rajput89}:
in order to perform stochastic integration, we need to single out one
coordinate to be time and introduce a filtration based definition of
the integrator $\La$. For notational convenience, we will write $\La
(A)$ instead of $\La(\Om\times A)$ in the following. As shown in
\citep{Rajput89}, Proposition 2.1 and Lemma~2.3, $\La$ induces a
characteristic triplet $(B,C,\nu)$ w.r.t. some truncation function
$\tau$ via the L\'evy--Khintchine formula:
\[
\bbe\bigl[\ee^{\ii u\La(A)}\bigr] = \exp \biggl(\ii u B(A) - \frac{u^2}{2}
C(A) + \int_\bbr\bigl(\ee^{\ii u y} - 1 - \ii u\tau(y)
\bigr) \nu(A,\dd y) \biggr),\qquad  A\in\cals,u\in\bbr.
\]
It is natural to ask how this notion of characteristics compares with
Theorem~\ref{candec}. Obviously, $\La$~is an orthogonal random
measure. In order that $\La$ has different times of discontinuity, it
suffices by independence to assume that $\La$ has no fixed times of
discontinuity. In this case, recalling the construction in the proof of
Theorem~\ref{candec} and using \citep{Sato04}, Theorem 3.2,  together with
\citep{Jacod03}, Theorem II.4.15, one readily sees that the two different
definitions of characteristics agree in the natural filtration of $\La
$. In particular, the canonical decomposition of $\La$ determines its
L\'evy--It\^o decomposition as derived in \citep{Pedersen03}.\looseness=-1

Consequently, the integrability criteria obtained in Theorem~\ref
{main} extend the corresponding result of \citep{Rajput89}, Theorem 2.7,
for deterministic functions (or, as used in \citep{BN11-2}, for
integrands which are independent of $\La$) to allow for predictable
integrands.
\end{Remark}

%
\begin{Remark}[(Martingale measures \citep{Walsh86})]\label{Walsh} In
\citep
{Walsh86}, a stochastic integration theory for predictable integrands
is developed with so-called worthy martingale measures as integrators.
The concept of worthiness is needed since a martingale measure in
Walsh's sense does not guarantee that it is a random measure in the
sense of Definition~\ref{ranmeas}. What is missing is, loosely
speaking, a joint \mbox{$\si$-additivity} condition in space and time; see
also the example in \citep{Walsh86}, page 305ff. The worthiness of a
martingale measure, that is, the existence of a dominating ($\si
$-additive) measure, turns it into a random measure.

In essence, the integration theory presented in \citep{Walsh86} for
worthy martingale measures is an $L^2$-theory similar to \citep
{Doob53,Ito44}, where the extension from simple to general integrands
is governed by a dominating measure. The latter also determines whether
a predictable function is integrable or not in terms of a
square-integrability condition; see \citep{Walsh86}, page 292. We see
the main advantages of the $L^2$-theory as follows: it does not require
the martingale measure to have different times of discontinuity, works
with fairly easy integrability conditions and produces stochastic
integrals again belonging to $L^2$. However, many interesting
integrators (e.g., stable noises) are not $L^2$-random measures.
Moreover, even if the integrator $M$ is an $L^2$-random measure, the
class $L^0(M)$ is usually considerably larger than the class $L^2(M)$.
Thus, in comparison to \citep{Walsh86}, it is the compensation of
these two shortages of the $L^2$-theory that constitutes the main
advantage of our
integrability conditions in Theorem~\ref{main}. We will come back to
this point in Section~\ref{s52}, where it is shown that in the study
of stochastic PDEs, solutions often do not exist in the $L^2$-sense but
in the $L^0$-sense.
\end{Remark}

%
\begin{Remark}[((Compensated) strict random measures \citep{Jacod03})]
Chapters I and II of \citep{Jacod03} are an established reference for
integration theory w.r.t. semimartingales. Moreover, they also cover
the integration theory w.r.t. strict random measures or compensated
strict random measures as follows: if $M$ is a strict random measure,
they define stochastic integrals w.r.t. $M$ path-by-path. More
precisely, a measurable function $H\dvtx \tilde\Om\to\bbr$ is
pathwise integrable w.r.t. $M$ if for a.e. ${\omega}\in\Om$
%
\begin{equation}
\label{strict} \int_{\bbr\times E} |H|({\omega},t,x) |M|({\omega}, \dd
t,\dd x) < \infty.
\end{equation}
If $\tilde M:=M-M^\mathrm{p}$ is the compensation of an integer-valued
strict random measure $M$, we have the following situation: let $H\in
\tilde\calp$ and introduce an auxiliary process by
%
\begin{equation}
\label{aux} \tilde H_t({\omega}):=\int_E
H({\omega},t,x) \tilde M\bigl({\omega},\{t\}\times\dd x\bigr), \qquad ({\omega},t)\in
\bar\Om,
\end{equation}
hereby setting $\tilde H_t({\omega}):=+\infty$ whenever \eqref{aux}
diverges. Then $H$ is integrable in the sense of \citep{Jacod03}, Definition
II.1.27,  if there exists a sequence of stopping times
$(T_n)_{n\in\bbn}$ with $T_n\uparrow+\infty$ a.s. and
%
\begin{equation}
\label{condcomp} \bbe \biggl[ \biggl(\sum_{-T_n\leq t
\leq T_n} \tilde
H_t^2 \biggr)^{1/2} \biggr] < \infty.
\end{equation}

How do these integrability conditions compare to those of Theorem~\ref
{main}? Obviously, pathwise integrability w.r.t. $M$ does not require
the integrand to be predictable. Furthermore, if $H$ is predictable and
\eqref{strict} holds, then the pathwise integral coincides with the
stochastic integral $H\cdot M$. Still, Theorem~\ref{main} provides a
useful extension in some situations: first, there are examples $H\in
L^0(M)$ which fail the condition \eqref{strict} (see the example at
the end of Remark~\ref{h0}). And second, given some specific $H$, it
may be difficult in general to determine whether \eqref{strict} holds
or not (e.g., if $M$ has no finite first moment). The characteristic
triplet that is used in Theorem~\ref{main} is often easier to handle
than $|M|$.

As for $\tilde M$ we have following situation: first, one should notice
that \eqref{condcomp} ensures integrability on \emph{finite}
intervals, whereas Theorem~\ref{main} is concerned with \emph{global}
integrability on $\bbr$. Second, even on finite intervals, the
conditions of Theorem~\ref{main} are more general than \eqref
{condcomp}, see \citep{Bichteler83}, Proposition 3.10. Finally, whereas
\eqref{condcomp} involves a localizing sequence of stopping times and
moment considerations, Theorem~\ref{main} relates integrability only
to the integrand itself and the characteristics of $\tilde M$, which is
often more convenient.
\end{Remark}

In order to illustrate condition \eqref{cond2} in Theorem~\ref{main},
we now discuss the example of a Gaussian random measure, which is white
in time but coloured in space. Such random measures are often
encountered as the driving noise of stochastic PDEs, see \citep
{Dalang99} and references therein.

\begin{Example}\label{LambdaH} Let $(M(\Om\times B))_{B\in\calb_\bb
(\bbr
^{1+d})}$ be a mean-zero Gaussian process whose covariance functional
for $B,B^\prime\in\calb_\bb(\bbr^{1+d})$ is given by
%
\begin{equation}
\label{CBB} C\bigl(B;B^\prime\bigr):=\bbe\bigl[M(\Om\times B)M\bigl(
\Om \times B^\prime\bigr)\bigr] = \int_\bbr\int
_{B(t)\times B^\prime(t)} f\bigl(x-x^\prime\bigr) \,\dd\bigl(x,
x^\prime\bigr) \,\dd t,
\end{equation}
where $B(t):=\{x\in\bbr^d\colon(t,x)\in B\}$.
For the existence of such a process, it is well known (\citep{Doob53}, Theorem
II.3.1), that $f\dvtx \bbr^d\to[0,\infty)$ must be a
symmetric and nonnegative definite function for which the integral on
the right-hand side of \eqref{CBB} exists. Under these conditions, $C$
defines a deterministic bimeasure which is symmetric in $B,B^\prime\in
\calb_\bb(\bbr^{1+d})$.

For the further procedure let $(\calf_t)_{t\in\bbr}$ be the natural
filtration of $M$ and set
\[
M\bigl(F\times(s,t]\times U\bigr) := 1_F M\bigl(\Om\times(s,t]\times U\bigr), \qquad F\in\calf
_s.
\]
By \citep{Bichteler83}, Theorem 2.25, $M$ can be extended to a random
measure on $\bbr\times\bbr^d$ provided that
\[
S_n \to0\qquad \mbox{pointwise},\qquad  |S_n|\leq|S|\quad
\Longrightarrow\quad \int S_n\, \dd M \to 0 \qquad \mbox{in } L^0
\]
for all step functions $S_n$ and $S$ over sets of the form $F\times
(s,t]\times U$ with $F\in\calf_s$, $s<t$ and $U\in\calb_\bb(\bbr
^d)$. Indeed, using obvious notation and observing that $1_F$ is
independent of $M(\Om\times(s,t]\times U)$ for $F\in\calf_s$ since
$M$ is white in time, we have
\begin{eqnarray*}
&& \bbe \biggl[ \biggl(\int S_n \,\dd M \biggr)^2 \biggr]\\
&&\quad = \sum_{i,j=1}^{r_n} a_i^n
a_j^n\bbe\bigl[M\bigl(A_i^n
\bigr)M\bigl(A_j^n\bigr)\bigr]
\\
&&\quad = \sum_{i,j=1}^{r_n} a_i^n
a_j^n P\bigl[F_i^n\bigr]P
\bigl[F_j^n\bigr] \Leb \bigl((s_i^n,t_i^n]
\cap(s_j^n,t_j^n]\bigr)
\int_{U_i^n\times U_j^n} f\bigl(x-x^\prime\bigr) \dd\bigl(x,
x^\prime\bigr)
\\
&&\quad = \int_{(\bbr^{1+d};\bbr^{1+d})} (\tilde S_n,\tilde
S_n) \,\dd C \to 0
\end{eqnarray*}
by dominated convergence \citep{Chang83}, Corollary 2.9. Here $\tilde S_n$
arises from $S_n$ by replacing $a_i^n$ with $a_i^n P[F_i^n]$.

Having established that $M$ is a random measure, let us derive its
characteristics. Obviously, $B$~and $\nu$ are identically $0$. It is
also easy to see that $C$ is the second characteristic of $M$: it is
clear for sets of the form $(s,t]\times U$, and extends to general sets
in $\calb_\bb(\bbr^{1+d})$ by dominated convergence. Therefore, as
shown in the proof of Theorem~\ref{charint}, $L^{1,0}(M)$ consists of
those $H\in\tilde\calp$ such that $(H;H)$ is strictly
$C$-integrable, or, equivalently,
%
\begin{equation}
\label{Cexample} \int_\bbr\int_{\bbr^d\times\bbr^d}
|H|(t,x)|H|\bigl(t,x^\prime\bigr)f\bigl(x-x^\prime\bigr) \,\dd
\bigl(x,x^\prime\bigr) \,\dd t < \infty \qquad \mbox{a.s.}
\end{equation}
The class $L^0(M)$, however, is the set of all $H\in\tilde\calp$
such that a.s. the inner integral in \eqref{Cexample} is finite for
a.e. $t\in\bbr$, and
%
\begin{equation}
\label{Cexample2} \int_\bbr\int_{\bbr^d\times\bbr^d}
H(t,x)H\bigl(t,x^\prime\bigr)f\bigl(x-x^\prime\bigr) \,\dd
\bigl(x,x^\prime\bigr) \,\dd t < \infty \qquad \mbox{a.s.}
\end{equation}
A (deterministic) function $H\in L^0(M)$ which is not in $L^{1,0}(M)$
is, for instance, given by $H(t,x):= th(x)$ where $h$ is chosen such
that
\[
\int_{\bbr^d\times\bbr^d} h(x)h\bigl(x^\prime\bigr)f
\bigl(x-x^\prime\bigr) \,\dd \bigl(x,x^\prime\bigr) = 0.
\]

One important example is a fractional correlation structure in space.
In this case, we have $f(x_1,\ldots,x_d)=\prod_{i=1}^d
|x_i|^{2H_i-2}$, where $H_i\in(1/2,1)$ is the Hurst index of the $i$th
coordinate. Then $L^0(M)$ can be interpreted as the extension of the
class $|\La_H|$ studied in \citep{Pipiras00} to several parameters
and stochastic integrands. However, in \citep{Pipiras00} as well as in
\citep{Basse12}, stochastic integrals are constructed for even larger
classes of integrands. These classes, denoted $\La_H$ or $\La_X$,
respectively, are obtained as limits of simple functions under
$L^2$-norms ($\|\cdot\|_{\La_H}$ and $\|\cdot\|_{\La_X}$,
resp.), which are defined via fractional derivatives or Fourier
transforms. In particular, the stochastic integrals defined via these
norms are no longer of It\^o type, that is, no dominated convergence
theorem holds for these stochastic integrands. Indeed, $L^{1,0}(M)$ is
the largest class of predictable integrands for which a dominated
convergence theorem
holds (see Theorem~\ref{L1p-char}), and $L^0(M)$ is its improper
extension to functions for which $H\cdot M$ is a finite measure.
\end{Example}

The investigation of multi-dimensional stochastic processes often
involves stochastic integrals where the integrand $H$ is a
matrix-valued predictable function and the integrator $M=(M^1,\ldots
,M^d)$ is a $d$-dimensional random measure, that is, $M^1, \ldots,
M^d$ are all random measures in the sense of Definition~\ref{ranmeas}
w.r.t. the same underlying filtration and the same sequence $(\tilde
O_k)_{k\in\bbn}$.
By considering each row of $H$ separately, we can assume for the
following that $H$ is an $\bbr^d$-valued predictable function. It is
obvious that the construction of stochastic integrals requires no more
techniques than those presented in Section~\ref{s2}. In fact,
replacing $E$ by $E^d$ reduces\linebreak[4]  the multivariate case to the univariate one.
However, there is a difference when we want to apply the canonical
decomposition as in Theorem~\ref{candec} or the integrability
conditions in Theorem~\ref{main}: in the multi-dimensional case, it is
not reasonable to assume that $M^i$ and $M^j$ for $i\neq j$ have
different times of discontinuity. Instead, one would define
$d$-dimensional characteristics $(B,C,\nu)$ for $M$, similar to \citep{Jacod03},
Chapter II,  or \citep{Basse13}, Theorem 3.1, and use these to
characterize integrability.

In the next theorem, we rephrase \ref{main} for the multivariate
setting. Since no novel arguments are needed, we omit its proof. We
will use the product notation in a self-explanatory way: for instance,
if $x, y \in\bbr^d$, $xy$ denotes their inner product; for $A\in
\tilde\calp_M$, $1_A\cdot M$ denotes the $d$-dimensional
semimartingale $(1_A\cdot M^1,\ldots,1_A\cdot M^d)$; $H\cdot M$
denotes $\sum_{i=1}^d H^i\cdot M^i$ for $H\in L^{1,0}(M)$ and is
suitably extended to $H\in L^0(M)$, cf. Section~\ref{s2}. Similarly,
given a matrix $\beta=(\beta^{ij})_{i,j=1}^d$ of bimeasures from
$\calf_1\times\calf_2\to\bbr$ and $\calf_i$-measurable functions
$f_i=(f^1_i,\ldots,f^d_i)$ for $i=1,2$, we define
\[
\int_{(A;B)} (f_1;f_2) \,\dd\beta:=\sum
_{i,j=1}^d \int_{(A;B)}
\bigl(f_1^i;f_2^j\bigr) \,\dd
\beta^{ij}, \qquad A\in\calf_1,B\in\calf_2,
\]
whenever the right-hand side exists.

Assume that $M$ has different times of discontinuity, which means that
$1_{O_k\times U_i} \cdot M$, $i=1, 2$, a.s. never jump at the same time
for all disjoint sets $U_1, U_2\in\cale_M$ and $k\in\bbn$. Given a
truncation function $\tau\dvtx \bbr^d\to\bbr^d$, define for
$A,B\in\tilde\calp_M$ and $V\in\calb_0(\bbr^d)$
%
\begin{eqnarray}\label{charmultvers}
\hspace*{-20pt}B(A)&:=&\mathfrak{B}(1_A\cdot M)_\infty, \qquad \mu(A,V):=
\mu^{1_A\cdot
M}(\bbr,V),\qquad  \nu(A,V):=\nu^{1_A\cdot M}(\bbr,V)
\nonumber
\\[-8pt]\\[-8pt]
\hspace*{-20pt}M^\cc(A)&:=&(1_A\cdot M)^\cc,\qquad
C^{ij}(A;B):=\bigl[\bigl(1_A\cdot M^i
\bigr)^\cc ,\bigl(1_B\cdot M^j
\bigr)^\cc\bigr]_\infty. \nonumber
\end{eqnarray}
As in Theorem~\ref{candec} $(B,C,\nu)$ can be extended to predictable
strict random (bi-)measures and give rise to the following canonical
decomposition of $M$:
%
\begin{eqnarray}
M(A)&=&B(A)+M^{\cc}(A) + \int_{\bbr\times E \times\bbr
^d}
1_A(t,x) \bigl(y-\tau(y)\bigr) \mu(\dd t,\dd x,\dd y)
\nonumber
\\[-8pt]\\[-8pt]
&&{}+ \int_{\bbr\times E \times\bbr^d} 1_A(t,x)\tau(y) (\mu-\nu ) (\dd
t,\dd x,\dd y),\qquad  A\in\tilde\calp_M.\nonumber
\end{eqnarray}
Moreover, there exist a positive predictable strict random measure
$A({\omega},\dd t,\dd x)$, a $\tilde\calp$-measurable $\bbr
^d$-valued function $b({\omega},t,x)$ and a transition kernel
$K({\omega},t,x,\dd y)$ from $(\tilde\Om,\tilde\calp)$ to $(\bbr
^d,\calb(\bbr^d))$ such that for all ${\omega}\in\Om$,
\[
B({\omega},\dd t,\dd x)=b({\omega},t,x) A({\omega},\dd t,\dd x), \qquad \nu ({\omega},
\dd t,\dd x, \dd y)=K({\omega},t,x,\dd y) A({\omega}, \dd t,\dd x).
\]

The multi-dimensional version of Theorem~\ref{main} reads as follows.

\begin{Theorem}\label{main-mult} Let $M$ be a $d$-dimensional random measure
with different times of discontinuity and $H\dvtx \tilde\Om\to\bbr
^d$ be a predictable function such that there exists a strictly
positive predictable process $K\dvtx \bar\Om\to\bbr$ with $HK\in
L^{1,0}(M)$.
Then $H\in L^0(M)$ if and only if each of the following conditions is
satisfied a.s.:
\begin{eqnarray*}
\int_{\bbr\times E} \biggl\llvert H(t,x)b(t,x) + \int
_{\bbr^d} \bigl[\tau\bigl(H(t,x)y\bigr)-H(t,x)\tau(y) \bigr]
K(t,x,\dd y)\biggr\rrvert A(\dd t,\dd x)& <& \infty,
\\
\int_\bbr K^{-2}_t \,\dd \biggl(\int
_{((-\infty,t]\times E;(-\infty
,t]\times E)} (HK;HK) \,\dd C \biggr)&<&\infty,
\\
\int_{\bbr\times E} \int_{\bbr^d} \bigl(1
\wedge\bigl|H(t,x)y\bigr|^2 \bigr) K(t,x,\dd y) A(\dd t,\dd x) &<& \infty.
\end{eqnarray*}
\end{Theorem}

\section{Ambit processes}\label{s5}

In this section, we present various applications, where the
integrability conditions of Theorem~\ref{main} are needed. Given a
filtered probability space satisfying the usual assumptions, our
examples are processes of the following form:
%
\begin{equation}
\label{ambit} Y(t,x):=\int_{\bbr\times\bbr^d} h(t,s;x,y) M(\dd s,\dd y), \qquad t
\in \bbr,x\in\bbr^d,
\end{equation}
where $h\dvtx \bbr\times\bbr\times\bbr^d\times\bbr^d \to\bbr$
is a deterministic measurable function and $M$ a random measure with
different times of discontinuity such that the integral in \eqref
{ambit} exists in the sense of \eqref{defint}. If the characteristics
of $M$ in the sense of Theorem~\ref{candec} are known, \eqref{ambit}
exists if and only if the conditions of Theorem~\ref{main} are
satisfied for each pair $(t,x)\in\bbr\times\bbr^d$. We call
processes of the form \eqref{ambit} \textit{ambit processes} although
the original definition in \citep{BN11-2} requires the random measure
to be a volatility modulated L\'evy basis, that is, $M=\si.\La$ where
$\La$ is a L\'evy basis and $\si\in\tilde\calp$. As already
explained in the \hyperref[s1]{Introduction}, this class of models is relevant in many
different areas of applications. In the following subsections, we
discuss two applications where interesting choices for $h$ and $M$ will
be presented.

\subsection{Stochastic PDEs}\label{s51}
The connection between stochastic PDEs and ambit processes is
exemplified in \citep{BN11-2} relying on the integration theory of
\citep{Rajput89} or \citep{Walsh86}. Let $U$ be an open subset of
$\bbr\times\bbr^d$ with boundary $\partial U$, $P$ a polynomial in
$1+d$ variables with constant coefficients and $M$ a random measure
with different times of discontinuity. The goal is to find a solution
$Z$ to the stochastic PDE
%
\begin{equation}
\label{spde} P(\partial_t,\partial_1,\ldots ,\partial
_d)Z(t,x)=\partial_t\,\partial_1\cdots
\partial_d M(t,x), \qquad (t,x)\in U,
\end{equation}
subjected to some boundary conditions on $\partial U$, where $\partial
_t\,\partial_1\cdots\partial_d M(t,x)$ is the formal derivative of
$M$, its noise. We want to apply the method of Green's function to our
random setting: first, we find a solution $Y$ to \eqref{spde} with
vanishing boundary conditions, then we find a solution $Y^\prime$ to
the homogeneous version of \eqref{spde} which satisfies the prescribed
boundary conditions, and finally we obtain a solution $Z$ by the sum of
$Y$ and $Y^\prime$. Since the problem of finding $Y^\prime$ is the
same as in ordinary PDE theory, we concentrate on finding $Y$. However,
since the noise of $M$ does not exist formally, there exists no
solution $Y^\prime$ in the strong sense. One standard approach based
on \citep{Walsh86}, Section~3,  is to interprete \eqref{spde} in weak
form and to define
%
\begin{equation}
\label{Ytx} Y(t,x):=\int_U G(t,s;x,y) M(\dd s,\dd y),\qquad
(t,x)\in U,
\end{equation}
as a solution, where $G$ is the Green's function for $P$ in the domain
$U$. Obviously, $Y$ is then an ambit process, where the integrand is
determined by the partial differential operator and the domain, and the
integrator is the driving noise of the stochastic PDE. Therefore,
Theorem~\ref{main} provides necessary and sufficient conditions for
the existence of $Y$. Let us stress again that, in contrast to \citep
{Rajput89} and \citep{Walsh86}, we need no distributional assumptions
on~$M$.

Finally, we want to come back to Remark~\ref{Walsh} and explain why
the $L^2$-approach is too stringent for stochastic PDEs. To this end,
we consider the stochastic heat equation in~$\bbr^d$:

\begin{Example}We take $P(t,x)=t-\sum_{i=1}^d x_i$, $U=(0,\infty
)\times\bbr
^d$ and $M=\si.\La$ where $\si$ is a predictable function and $\La$
a L\'evy basis with characteristics $(0,\Sigma\,\dd t \,\dd x, \nu
(\dd\xi) \,\dd t\,\dd x)$, where $\Sigma\geq0$ and $\nu$ is a symmetric L\'
evy measure. \citep{Walsh86}, Section~3,  considers a similar equation
with $\nu=0$. The Green's function for $P$ and $U$ is the heat kernel
\[
G(t,s;x,y)=\frac{\exp(-|x-y|^2/(4(t-s)))}{(4\uppi(t-s))^{d/2}}1_{\{
0<s< t\}}, \qquad s,t>0,x,y\in\bbr^d.
\]
Since for all $(t,x)\in U$ the kernel $G(t,\cdot;x,\cdot)\in L^p(U)$
if and only if $p<1+2/d$, it is square-integrable only for $d=1$.
Therefore, in the $L^2$-approach function-valued solutions only exist
for $d=1$. However, if $\Sigma=0$, a sufficient condition for \eqref
{cond3} and thus the existence of \eqref{Ytx} is
%
\begin{eqnarray}
\label{heatex} \int_0^t \int
_{\bbr^d} \bigl|G(t,s;x,y)\si(s,y)\bigr|^p \,\dd s \,\dd y& < &\infty\qquad
\mbox{a.s.},\quad  (t,x)\in U,\nonumber\\[-8pt]\\[-8pt]
 \int_{[-1,1]} |\xi|^p \nu(
\dd\xi)& < &\infty\nonumber
\end{eqnarray}
for some $p\in[0,2)$. For instance, if $\si$ is stationary in $U$
with finite $p$th moment, \eqref{heatex} becomes
\[
\int_{[-1,1]} |\xi|^p \nu(\dd\xi)<\infty\qquad \mbox{for
some } p<1+2/d.
\]
In particular, we see that function-valued solutions exist in arbitrary
dimensions, which cannot be ``detected'' in the $L^2$-framework, even
for integrators which are $L^2$-random measures.
\end{Example}

The stochastic heat equation or similar equations driven by
non-Gaussian noise have already been studied in a series of papers,
for example, \citep{Albeverio98,Applebaum00,Mueller98,Mytnik02,SLB98}, partly
also extending Walsh's approach beyond the $L^2$-framework. Although
they do not only consider the linear case \eqref{spde}, there are
always limitations in dimension (e.g., only $d=1$) or noise type (e.g.,
only stable noise without volatility modulation). Thus, in the linear
case, Theorem~\ref{main} provides a unifying extension of the
corresponding results in the given references.

\subsection{Superposition of stochastic volatility models}\label{s52}

In this subsection, we give examples of ambit processes, where the
spatial component in the stochastic integral has the meaning of a
parameter space. First, we discuss one possibility of constructing a
superposition of COGARCH processes, following \citep{Behme13}. The
COGARCH model of \citep{Klueppelberg04} itself is designed as a
continuous-time version of the celebrated GARCH model and is defined as follows:
Let $(L_t)_{t\in\bbr}$ be a two-sided L\'evy process with L\'evy
measure $\nu_L$. Given $\beta,\eta>0$ the COGARCH model $(V^\vp
,G^\vp)$ with parameter $\vp\geq0$ is given by the equations
%
\begin{eqnarray}
\label
{cog-price}\dd G^\vp_t &=& \sqrt{V^\vp_{t-}}
\,\dd L_t, \qquad G^\vp_0=0,
\\
\label{cog-sde}\dd V^\vp_t &=& \bigl(\beta-\eta V^\vp_t
\bigr) \,\dd t + \vp V^\vp_{t-} \,\dd S_t,\qquad  t\in
\bbr,
\end{eqnarray}
where $S:=[L]^\dd$ denotes the pure-jump part of the quadratic
variation of $L$. By \citep{Klueppelberg04}, Theorem~3.1, \eqref
{cog-sde} has a stationary solution if and only if
%
\begin{equation}
\label{condphi} \int_{\bbr_+} \log\bigl(1+\vp y^2
\bigr) \nu _L(\dd y) <\eta.
\end{equation}
Let us denote the collection of all $\vp\geq0$ satisfying \eqref
{condphi} by $\Phi$, which by \eqref{condphi} must be of the form
$[0,\vp_{\mathrm{max}})$ with some $0<\vp_{\mathrm{max}}<\infty$.
Although the COGARCH model essentially reproduces the same stylized
features as the GARCH model, there are two unsatisfactory aspects:
\begin{enumerate}[(2)]
\item[(1)] Right from the definition, the COGARCH shows a deterministic
relationship between volatility and price jumps, an effect shared by
many continuous-time stochastic volatility models \citep{Jacod12}.
More precisely, we have
%
\begin{equation}
\label{fixrel} \Delta V_t^\vp= \vp V_{t-}^\vp(
\Delta L_t)^2 = \vp\bigl(\Delta G^\vp_t
\bigr)^2,\qquad  t\in\bbr.
\end{equation}
A realistic stochastic volatility model should allow for different
scale parameters $\vp$.
\item[(2)] The autocovariance function of the COGARCH volatility is, when
existent and $\vp>0$, always of exponential type: $\cov[V^\vp
_t,V^\vp_{t+h}]= C\ee^{-ah}$ for $h\geq0$, $t\in\bbr$ and some
constants $C,a>0$. A more flexible autocovariance structure is desirable.
\end{enumerate}

In \citep{Behme13}, three approaches to construct superpositions of
COGARCH processes (supCO\-GARCH) with different values of $\vp$ are
suggested in order to obtain a stochastic volatility model keeping the
desirable features of the COGARCH but avoiding the two disadvantages
mentioned above. One of them is the following: With $\beta$ and $\eta
$ remaining constant, take a L\'evy basis $\La$ on $\bbr\times\Phi$
with characteristics $(b \dd t \pi(\dd\vp), \Sigma\dd t \pi(\dd
\vp), \nu_L(\dd y)\,\dd t \pi(\dd\vp))$, where $b\in\bbr$, $\Sigma
\geq0$, $\pi$ is a probability measure on $\Phi$ and $\nu_L$ the L\'
evy measure of the L\'evy process given by
\[
L_t:=\La^L((0,t]\times\Phi), \qquad t\geq0,\qquad
L_t:=-\La^L((-t,0]\times\Phi ),\qquad  t<0,
\]
Furthermore, define another L\'evy basis by $\La^S(\dd t,\dd\vp):=
\int_\bbr y^2 \mu^\La(\dd t,\dd\vp,\dd y)$,
where $\mu^\La$ is the jump measure of $\La$ as in Theorem~\ref
{candec}. Next define $V^\vp$ for each $\vp\in\Phi$ as the COGARCH
volatility process driven by $L$ with parameter $\vp$. Motivated by
\eqref{cog-sde}, the \textit{supCOGARCH} $\bar V$ is now defined by the
stochastic differential equation
%
\begin{equation}
\label{supcog-sde} \dd\bar V_t = (\beta- \eta\bar V_t) \dd
t + \int_{\Phi} \vp V^\vp_{t-} \La(\dd
t, \dd\vp),\qquad  t\in\bbr.
\end{equation}
As shown in \citep{Behme13}, Proposition 3.15, \eqref{supcog-sde} has a
unique solution given by
%
\begin{equation}
\label{supcog} \bar V_t=\frac{\beta}{\eta}+\int_{-\infty}^t
\int_\Phi\ee^{-\eta(t-s)} \vp V^\vp_{s-}
\La(\dd s,\dd\vp), \qquad t\in\bbr,
\end{equation}
such that $\bar V$ is an ambit process as in \eqref{ambit}.

Here the integrability conditions of Section~\ref{s4} come into play.
Immediately from Theorem~\ref{main} and Remark~\ref{h0}, we obtain
the following corollary.

\begin{Corollary}
The supCOGARCH $\bar V$ as in \eqref{supcog} exists if and only if
%
\begin{equation}
\label{ex-cog} \int_{\bbr_+} \int_\Phi\int
_{\bbr
_+} 1\wedge \bigl(y^2 \vp\ee^{-\eta s}
V_s^\vp \bigr) \nu_L(\dd y) \pi(\dd\vp ) \dd
s<\infty\qquad \mbox{a.s.}
\end{equation}
\end{Corollary}

In particular, the supCOGARCH \eqref{supcog} provides an example where
the stochastic volatility process $\si(s,\vp):=\vp V^\vp_{s-}$ is
\emph{not} independent of the underlying L\'evy basis $\La$. So the
conditions of \citep{Rajput89}, Theorem 2.7,  are not applicable. For
further properties of the supCOGARCH, in particular regarding its jump
behaviour, autocovariance structure etc., we refer to \citep{Behme13}.

Finally, let us comment on superpositions of other stochastic
volatility models.

%
\begin{Remark}
The usage of Ornstein--Uhlenbeck processes in stochastic volatility
modelling has become popular through the Barndorff-Nielsen--Shephard
model \citep{BN01a}. A natural extension is given by the CARMA
stochastic volatility model \citep{Todorov06}, which generates a more
flexible autocovariance structure. Another generalization of the BNS
model is obtained via a superposition of OU processes with different
memory parameters leading to the class of supOU processes \citep
{BN01}. This method does not only yield a more general second-order
structure but can also generate long-memory processes; cf. \citep
{BN01,Fasen07}. A similar technique was used in \citep
{BN11-3,Marquardt07} to construct supCARMA processes, again leading to
a possible long-range dependent process.

Note that in all these models the driving noise is assumed to have
stationary independent increments, which is certainly a model
restriction. Therefore, \citep{BN12} suggests a volatility modulation
of this noise to obtain a greater model flexibility. In this way, it is
possible to generate a volatility clustering effect, similar to the
behaviour of the (sup)COGARCH. Without volatility modulation, supOU or
supCARMA processes are defined as stochastic integrals of deterministic
kernel functions w.r.t. a L\'evy basis, so the approach of \citep
{Rajput89} is sufficient. Theorem~\ref{main} now enables us to replace
$\La$ by a volatility modulated L\'evy basis $\si.\La$ with a
possible dependence structure between $\si$ and $\La$.
\end{Remark}


\section*{Acknowledgements}
We take pleasure in thanking Jean Jacod for extremely useful
discussions and his advice on this topic. We also thank two referees
and an Associate Editor for their constructive comments which led to a
considerable improvement of the paper. The first author acknowledges
support from the graduate program TopMath of the Elite Network of
Bavaria and the TopMath Graduate Center of TUM Graduate School at
Technische Universit\"at M\"unchen.



\printhistory
\end{document}